\documentclass[12pt]{amsart}
\usepackage[utf8]{inputenc}
\usepackage{amsmath,amssymb,amsthm}
\usepackage{enumerate, euler, palatino}
\usepackage{mathrsfs}
\usepackage{tikz-cd}
\usepackage[colorlinks=true]{hyperref}
\hypersetup{urlcolor=blue, citecolor=red , linkcolor= blue}
\usepackage{cite}
\usepackage[a4paper,margin=2.5cm]{geometry}

\theoremstyle{plain}
\newtheorem{theorem}{Theorem}[section]
\newtheorem{proposition}[theorem]{Proposition}
\newtheorem{lemma}[theorem]{Lemma}
\newtheorem{corollary}[theorem]{Corollary}
\theoremstyle{definition}

\newtheorem{example}[theorem]{Example}
\newtheorem{remark}[theorem]{Remark}

\begin{document}

\title[Taylor Invertibility and the Completion Problem]{Taylor Invertibility and the Completion Problem for Commuting Tuples of Upper Triangular Operator Matrices}
\author[C. Benhida]{Chafiq Benhida}
\address{UFR de Math\'ematiques, Universit\'e des Sciences et Technologies de Lille, F-59655 Villeneuve-d'Ascq Cedex, France}
\email{\url{chafiq.benhida@univ-lille.fr}}

\author[H. Stankovi\'c]{Hranislav Stankovi\'c}
\address{Faculty of Electronic Engineering, University of Ni\v s, Aleksandra Medvedeva 4, 18000 Ni\v s, Serbia}
\email{\url{hranislav.stankovic@elfak.ni.ac.rs}}
\date{\today}

\begin{abstract}
We study commuting $n$-tuples of upper triangular operator matrices $\mathbf{M}_{\mathbf{C}} = \bigl(\begin{smallmatrix} \mathbf{A} & \mathbf{C} \\ \mathbf{0} & \mathbf{B} \end{smallmatrix}\bigr)$ acting on a direct sum of Banach spaces. Using the long exact sequence in cohomology induced by the Koszul complexes, we show that $\mathbf{M}_{\mathbf{C}}$ is Taylor invertible if and only if the connecting homomorphisms $H^k(\mathbf{B}) \to H^{k+1}(\mathbf{A})$ are linear bijections in every degree, and that under natural closed-range hypotheses on the diagonal tuples these bijections are automatically topological isomorphisms. We determine which closed-range properties of the diagonal tuples are forced by the Taylor invertibility of $\mathbf{M}_{\mathbf{C}}$ and show by an example that the remaining ones fail in general. We then study the multivariable invertible completion problem, that is, the existence of a coupling tuple $\mathbf{C}$ for which $\mathbf{M}_{\mathbf{C}}$ is Taylor invertible. We reduce it to the realization of a family of cohomological isomorphisms by a single coupling tuple, show that the cohomology of $\mathbf{B}$ and the shifted cohomology of $\mathbf{A}$ must be isomorphic as graded modules over the exterior algebra, and exhibit diagonal tuples whose cohomology spaces are topologically isomorphic in every degree but which admit no Taylor invertible completion, a phenomenon with no one-dimensional counterpart.
\end{abstract}

\maketitle

\bigskip
\section{Introduction}\label{sec:intro}
\bigskip
In recent years, an extensive study of $2 \times 2$ upper triangular operator matrices has been carried out. This research was motivated primarily by the following fact: if $A \in \mathfrak{B}(\mathcal{X})$ and $\mathcal{S}$ is a closed, complemented $A$-invariant subspace of $\mathcal{X}$, then $A$ may be expressed in the block form
\begin{equation}
    A = \begin{pmatrix} * & * \\ 0 & * \end{pmatrix} : \begin{pmatrix} \mathcal{S} \\ \mathcal{P} \end{pmatrix} \to \begin{pmatrix} \mathcal{S} \\ \mathcal{P} \end{pmatrix},
\end{equation}
where $\mathcal{X} = \mathcal{S} \oplus \mathcal{P}$, with $\mathcal{P}$ closed.
Consequently, many recently published papers have focused on the issues of invertibility and regularity for an upper triangular operator matrix
\begin{equation} \label{eq:MC}
    M_C = \begin{pmatrix} A & C \\ 0 & B \end{pmatrix} : \begin{pmatrix} \mathcal{X} \\ \mathcal{Y} \end{pmatrix} \to \begin{pmatrix} \mathcal{X} \\ \mathcal{Y} \end{pmatrix},
\end{equation}
as well as on various types of spectra of $M_C$. A particular problem originating from these investigations deals with completions of the partially specified operator matrix $\begin{pmatrix} A & ? \\ 0 & B \end{pmatrix}$ in order to obtain an operator with some requested property. For given $A \in \mathfrak{B}(\mathcal{X})$ and $B \in \mathfrak{B}(\mathcal{Y})$, the following three questions repeatedly arise:
\begin{itemize}
    \item \textbf{Question 1:} Is there an operator $C \in \mathfrak{B}(\mathcal{Y}, \mathcal{X})$ such that $M_C$ is of a certain type (invertible, right invertible, left invertible, regular, \dots)?
    \item \textbf{Question 2:} For given operators $A$ and $B$, can we compute $\bigcap_{C} \sigma_*(M_C)$, where the intersection is taken over all $C \in \mathfrak{B}(\mathcal{Y}, \mathcal{X})$ and $\sigma_*$ is any type of spectra (essential, left, right, Weyl, Browder, etc.)?
    \item \textbf{Question 3:} Is there an operator $C' \in \mathfrak{B}(\mathcal{Y}, \mathcal{X})$ such that $\sigma_*(M_{C'}) = \bigcap_{C} \sigma_*(M_C)$, the intersection being taken over all $C \in \mathfrak{B}(\mathcal{Y}, \mathcal{X})$?
\end{itemize}
In general, the methodology for answering these questions in the single-variable case consists of finding necessary and sufficient conditions for the existence of $C$, followed by a description of the set of all such operators. For a comprehensive study of these problems, see \cite{Cvetkovic-Ilic22}.
The central result for the invertible completion is summarized in the following theorem, which is given in \cite{Cvetkovic-Ilic15,DuPan94,HanLeeLee00}.
\begin{theorem}\label{thm:1D_completion}
    Let $A \in \mathfrak{B}(\mathcal{X})$ and $B \in \mathfrak{B}(\mathcal{Y})$ be given operators. The operator matrix $M_C$ is invertible for some $C \in \mathfrak{B}(\mathcal{Y}, \mathcal{X})$ if and only if:
    \begin{enumerate}[\rm($i$)]
        \item $A$ is left invertible.
        \item $B$ is right invertible.
        \item $\mathcal{N}(B) \cong \mathcal{X}/\mathcal{R}(A)$ (as Banach spaces).
    \end{enumerate}
\end{theorem}
\bigskip
While the single-operator completion problem is well understood, generalizing these concepts to the multivariable setting, that is, to commuting $n$-tuples of operators, presents significant algebraic and topological challenges. The appropriate notion of joint invertibility in this context is provided by the Taylor spectrum, which is characterized by the exactness of the Koszul complex.

The relation between the Taylor invertibility of an upper triangular tuple and that of its diagonal entries is classical. Taylor \cite{Taylor70} showed, by means of the long exact sequence in cohomology, that for a commuting tuple, its restriction to a closed invariant subspace and the induced quotient tuple, the union of any two of the three Taylor spectra contains the third. In the Hilbert space setting, Curto \cite{Curto80} demonstrated that if an $n$-tuple of operators $\mathbf{A} = (A_1, \dots, A_n)$ and an $n$-tuple $\mathbf{B} = (B_1, \dots, B_n)$ are Taylor invertible, and if the operator matrix $\mathbf{M}_{\mathbf{C}} = \begin{pmatrix} \mathbf{A} & \mathbf{C} \\ \mathbf{0} & \mathbf{B} \end{pmatrix}$ forms a commuting $n$-tuple, then $\mathbf{M}_{\mathbf{C}}$ itself is Taylor invertible, and that if $\mathbf{A}$ and $\mathbf{M}_{\mathbf{C}}$ are Taylor invertible, then so is $\mathbf{B}$ (see also \cite{Curto81} for a systematic study of invertible and Fredholm $n$-tuples by means of Koszul complexes).

A point which runs through the whole paper is the distinction between algebraic and topological properties of Koszul cohomology. The Taylor invertibility of $\mathbf{M}_{\mathbf{C}}$ forces the boundary maps of its Koszul complex to have closed ranges, but this property need not be inherited by the diagonal tuples: we give a Taylor invertible pair $\mathbf{M}_{\mathbf{C}}$ for which some of the cohomology spaces of $\mathbf{A}$ and $\mathbf{B}$ fail to be Banach spaces. Accordingly, the long exact sequence in cohomology is used in the category of vector spaces, while closed-range hypotheses are introduced whenever the arguments require the Banach space structure of cohomology or concern the transfer of closed-range properties. This is a feature of the method rather than an obstruction to the completion problem, since a Taylor invertible completion may well exist in such a situation. A genuine obstruction appears at the algebraic level: the cohomology of $\mathbf{B}$ and the shifted cohomology of $\mathbf{A}$ must be isomorphic as graded modules over the exterior algebra, and we give diagonal tuples whose cohomology spaces are isomorphic, even topologically, in every degree, but whose module structures are not, so that no Taylor invertible completion exists.

Building upon this foundation, the present paper aims to bridge the gap between the single-variable completion problem and multivariable spectral theory. We investigate the conditions under which there exists an $n$-tuple $\mathbf{C} = (C_1, \dots, C_n)$ such that the commuting operator matrix $\mathbf{M}_{\mathbf{C}}$ is Taylor invertible. The homological viewpoint adopted here is a multivariable counterpart of the chain-complex approach to the $2 \times 2$ invertible completion problem developed in \cite{Stankovic1D}.

The paper is organized as follows. Section \ref{sec:prelim} collects the necessary background on Koszul complexes and Taylor invertibility. In Section \ref{sec:char} we determine which closed-range properties of the diagonal tuples are forced by the Taylor invertibility of $\mathbf{M}_{\mathbf{C}}$ (Theorem \ref{thm:AB_range_closed}), provide sufficient conditions for the closedness of the diagonal ranges in the remaining degrees (Theorems \ref{thm:B_range_closed} and \ref{thm:A_range_closed}), and show by an example that this closedness fails in general (Example \ref{ex:non_closed}). We then establish the long exact sequence associated with the triangular structure and characterize the Taylor invertibility of $\mathbf{M}_{\mathbf{C}}$ in terms of the connecting homomorphisms (Theorem \ref{thm:taylor_char}). As a corollary we recover, in the Banach space setting, the two-out-of-three principle of Taylor and Curto together with the accompanying spectral identity (Corollary \ref{cor:two_out_of_three}). Section \ref{sec:completion} is devoted to the multivariable invertible completion problem. Theorem \ref{thm:completion} reduces it to the realization of a family of cohomological isomorphisms by a single coupling tuple, Theorem \ref{thm:module} shows that the connecting homomorphisms are compatible with the natural action of the exterior algebra on Koszul cohomology, which yields a necessary condition for the existence of a completion in terms of graded modules, and Example \ref{ex:no_completion} exhibits diagonal tuples whose cohomology spaces are topologically isomorphic in every degree and which nevertheless admit no Taylor invertible completion.
\bigskip
\section{Preliminaries: Koszul Complexes and Taylor Invertibility}\label{sec:prelim}
\bigskip
Throughout the paper, $\mathcal{X}$ and $\mathcal{Y}$ denote complex Banach spaces, $\mathfrak{B}(\mathcal{X}, \mathcal{Y})$ the space of bounded linear operators from $\mathcal{X}$ to $\mathcal{Y}$, and $\mathfrak{B}(\mathcal{X}) := \mathfrak{B}(\mathcal{X}, \mathcal{X})$. For $T \in \mathfrak{B}(\mathcal{X}, \mathcal{Y})$, we write $\mathcal{N}(T) = \ker T$ and $\mathcal{R}(T)$ for the kernel and the range of $T$, respectively.

Let $\mathbf{e} = \{e_1, e_2, \dots, e_n\}$ be indeterminates and define $\Lambda_n[\mathbf{e}]$ to be the \textit{exterior algebra} on the generators $e_1, e_2, \dots, e_n$. This is the linear space over the complex plane $\mathbb{C}$ endowed with an anti-commutative exterior product $e_i \wedge e_j = -e_j \wedge e_i$ ($1 \leqslant i, j \leqslant n$). For $F = \{i_1, \dots, i_p\} \subseteq \{1, \dots, n\}$ with $i_1 < \dots < i_p$, we write $e_F = e_{i_1} \wedge \dots \wedge e_{i_p}$. The exterior algebra over $\mathbb{C}$ is then given by
\[
\Lambda_n[\mathbf{e}] = \left\{ \sum_{F} \alpha_F e_F : e_F = e_{i_1} \wedge \dots \wedge e_{i_p} \text{ and } \alpha_F \in \mathbb{C} \right\}.
\]
We let $e_{\emptyset}$ denote the identity element for the exterior product. If we denote $\Lambda_n^k[\mathbf{e}] = \left\{ \sum_{|F|=k} \alpha_F e_F : \alpha_F \in \mathbb{C} \right\}$, where $|F|$ is the cardinality of $F$, then clearly $\dim \Lambda_n^k[\mathbf{e}] = \binom{n}{k}$ for every $k \leqslant n$, $\Lambda_n^k[\mathbf{e}] \wedge \Lambda_n^l[\mathbf{e}] \subseteq \Lambda_n^{k+l}[\mathbf{e}]$ and $\Lambda_n[\mathbf{e}] = \bigoplus_{k=0}^n \Lambda_n^k[\mathbf{e}]$.
Given a Banach space $\mathcal{X}$, the exterior algebra over $\mathcal{X}$ is defined as:
\[
\Lambda_n[\mathbf{e}, \mathcal{X}] = \left\{ \sum_{F} x_F \otimes e_F : e_F = e_{i_1} \wedge \dots \wedge e_{i_p} \text{ and } x_F \in \mathcal{X} \right\}.
\]
The subspaces $\Lambda_n^p[\mathbf{e}, \mathcal{X}] = \left\{ \sum_{|F|=p} x_F \otimes e_F : x_F \in \mathcal{X} \right\}$, for $p \leqslant n$, are defined in a similar way. Naturally, $\Lambda_n^k[\mathbf{e}, \mathcal{X}]$ can be identified with $\mathcal{X}^{\binom{n}{k}}$, and we endow it with the norm
\[
\Bigl\| \sum_{|F|=k} x_F \otimes e_F \Bigr\| := \max_{|F|=k} \|x_F\|;
\]
any other choice of equivalent product norms would serve equally well. On a direct sum $\mathcal{X} \oplus \mathcal{Y}$ we use the norm $\|(x,y)\| := \max\{\|x\|, \|y\|\}$.
Since no confusion is possible, we will omit references to the indeterminates $e_1, e_2, \dots, e_n$ and write $\Lambda^k[\mathcal{X}]$ and $\Lambda_n[\mathcal{X}]$ for $\Lambda_n^k[\mathbf{e}, \mathcal{X}]$ and $\Lambda_n[\mathbf{e}, \mathcal{X}]$, respectively.
If $T \in \mathfrak{B}(\mathcal{X})$ (or, more generally, $T \in \mathfrak{B}(\mathcal{X}, \mathcal{Y})$), we use the same symbol $T$ to denote the operator defined on $\Lambda_n[\mathcal{X}]$ by
\[
T \left( \sum_{F} x_F \otimes e_F \right) = \sum_{F} T x_F \otimes e_F.
\]
For each $i \in \{1, 2, \dots, n\}$, let $E_i : \Lambda_n[\mathcal{X}] \to \Lambda_n[\mathcal{X}]$ be the left multiplication operator by $e_i$, defined by $E_i(x_F \otimes e_F) = x_F \otimes (e_i \wedge e_F)$. This is commonly called the \textit{creation operator}. With any commuting $n$-tuple $\mathbf{T} = (T_1, \dots, T_n)$, we associate the linear mapping defined over $\Lambda_n[\mathcal{X}]$ by
\[
\delta_{\mathbf{T}} = \sum_{i=1}^n T_i \otimes E_i : \sum_{F} x_F \otimes e_F \mapsto \sum_{F} \sum_{i=1}^n T_i x_F \otimes (e_i \wedge e_F).
\]
Clearly,
\begin{align*}
    \delta_{\mathbf{T}}^2(x_F \otimes e_F) &= \sum_{i,j=1}^n T_j T_i x_F \otimes (e_j \wedge e_i \wedge e_F) \\&= \sum_{i<j} T_iT_j x_F \otimes \left((e_j \wedge e_i +e_i \wedge e_j)\wedge e_F\right) = 0.
\end{align*}
It follows that $\delta_{\mathbf{T}}^2 = 0$ and hence $\mathcal{R}(\delta_{\mathbf{T}}) \subseteq \ker(\delta_{\mathbf{T}})$. Since $\delta_{\mathbf{T}}$ maps $\Lambda^{k-1}[\mathcal{X}]$ to $\Lambda^k[\mathcal{X}]$, the condition $\delta_{\mathbf{T}}^2 = 0$ is equivalent to $\mathcal{R}(\delta_{\mathbf{T}}|_{\Lambda^{k-1}[\mathcal{X}]}) \subseteq \ker(\delta_{\mathbf{T}}|_{\Lambda^k[\mathcal{X}]})$ for every $k \leqslant n$.
Set $\delta_{\mathbf{T}}^k := \delta_{\mathbf{T}}|_{\Lambda^k[\mathcal{X}]}$. We construct a cochain complex $\mathscr{K}(\mathbf{T})$, called the \textit{Koszul complex} associated with $\mathbf{T}$ on $\mathcal{X}$ as follows:
\[
\mathscr{K}(\mathbf{T}) : 0 \xrightarrow{\delta_{\mathbf{T}}^{-1}} \Lambda_n^0[\mathcal{X}] \xrightarrow{\delta_{\mathbf{T}}^0} \Lambda_n^1[\mathcal{X}] \xrightarrow{\delta_{\mathbf{T}}^1} \dots \xrightarrow{\delta_{\mathbf{T}}^{n-1}} \Lambda_n^n[\mathcal{X}] \xrightarrow{\delta_{\mathbf{T}}^n} 0.
\]
Following Taylor \cite{Taylor70}, the operator tuple $\mathbf{T}$ is said to be \textit{non-singular}, or \textit{Taylor invertible}, if $\ker \delta_{\mathbf{T}}^k = \mathcal{R}( \delta_{\mathbf{T}}^{k-1})$ for $k = 0, \dots, n$, which is equivalent to the statement that the associated Koszul complex is exact. The \textit{Taylor spectrum} of $\mathbf{T}$ on $\mathcal{X}$ is then defined as the set
\[
\sigma_T(\mathbf{T}) = \{ \lambda \in \mathbb{C}^n : \mathscr{K}(\mathbf{T} - \lambda) \text{ is not exact} \}.
\]
For the general theory of the Taylor spectrum and of invertible and Fredholm $n$-tuples, we refer the reader to \cite{Curto81,EschmeierPutinar96,Muller07}.
For a fixed $k$, we denote by
\[
H^k(\mathbf{T}) = \ker \delta_{\mathbf{T}}^k / \mathcal{R}(\delta_{\mathbf{T}}^{k-1})
\]
the $k$-th cohomology group of the Koszul complex $\mathscr{K}(\mathbf{T})$. Thus, $\mathbf{T}$ is Taylor invertible if and only if $H^k(\mathbf{T}) = \{0\}$ for all $k$; note that in this case each range $\mathcal{R}(\delta_{\mathbf{T}}^{k-1}) = \ker \delta_{\mathbf{T}}^{k}$ is automatically closed. In general, whenever the range $\mathcal{R}(\delta_{\mathbf{T}}^{k-1})$ is closed, the quotient $H^k(\mathbf{T})$ is a Banach space with respect to the quotient norm. If $\mathcal{R}(\delta_{\mathbf{T}}^{k-1})$ is not closed, the quotient topology on $H^k(\mathbf{T})$ is not Hausdorff, and $H^k(\mathbf{T})$ is regarded merely as a vector space. This algebraic reading is the default throughout the paper; the Banach space structure of a cohomology space is used only where the relevant range is known to be closed, and this will always be stated explicitly.

Throughout the paper, we adopt the convention that $\Lambda^k[\mathcal{X}] = \{0\}$, $\delta_{\mathbf{T}}^k = 0$, and hence $H^k(\mathbf{T}) = \{0\}$, for $k < 0$ and $k > n$. Accordingly, statements quantified over all $k \in \mathbb{Z}$ involve only finitely many nontrivial degrees, and the degenerate degrees carry precise content: for instance, $H^0(\mathbf{T}) = \ker \delta_{\mathbf{T}}^0$ and $H^n(\mathbf{T}) = \Lambda^n[\mathcal{X}]/\mathcal{R}(\delta_{\mathbf{T}}^{n-1})$.
\bigskip
\section{Characterization of Taylor Invertibility for Operator Matrix Tuples}\label{sec:char}
\bigskip
Let $\mathcal{X} \oplus \mathcal{Y}$ be a Banach space decomposition and consider an $n$-tuple of operators $\mathbf{M}=(M_1, \dots, M_n)$ acting on $\mathcal{X} \oplus \mathcal{Y}$, where each $M_i$ is an upper triangular operator matrix of the form
\begin{equation}\label{eq:M_tuple}
    M_i = \begin{pmatrix} A_i & C_i \\ 0 & B_i \end{pmatrix} \in \mathfrak{B}(\mathcal{X} \oplus \mathcal{Y}).
\end{equation}
The fundamental assumption that $\mathbf{M}$ is a commuting $n$-tuple is equivalent to the following statements:
\begin{itemize}
    \item $\mathbf{A}=(A_1, \dots, A_n)$ is a commuting $n$-tuple of bounded linear operators on $\mathcal{X}$;
    \item $\mathbf{B}=(B_1, \dots, B_n)$ is a commuting $n$-tuple of bounded linear operators on $\mathcal{Y}$;
    \item The coupling operators $\mathbf{C}=(C_1, \dots, C_n)$ must satisfy the following relation for all $i, j$:
    \begin{equation}\label{eq:coupling_relation}
        A_i C_j - A_j C_i = C_j B_i - C_i B_j.
    \end{equation}
\end{itemize}
In the language of Koszul complexes, the relation \eqref{eq:coupling_relation} is precisely the condition required for the boundary operator $\delta_{\mathbf{M}}$ to satisfy $\delta^2_{\mathbf{M}} = 0$.
First, for each $k \in \{0, \dots, n\}$, we naturally identify the space $\Lambda^k[\mathcal{X} \oplus \mathcal{Y}]$ with the direct sum $\Lambda^k[\mathcal{X}] \oplus \Lambda^k[\mathcal{Y}]$ via the isometric isomorphism
\begin{equation}
    \left( \sum_{|F|=k} \binom{x_F}{y_F} \otimes e_F \right) \cong \left( \sum_{|F|=k} x_F \otimes e_F, \sum_{|F|=k} y_F \otimes e_F \right)
\end{equation}
(with the norm conventions of Section \ref{sec:prelim}). Under this identification, an element $\zeta \in \Lambda^k[\mathcal{X} \oplus \mathcal{Y}]$ can be represented as a column vector $\begin{pmatrix} \xi \\ \eta \end{pmatrix}$, where $\xi \in \Lambda^k[\mathcal{X}]$ and $\eta \in \Lambda^k[\mathcal{Y}]$.
Representing the boundary operator $\delta_{\mathbf{M}}^k$ in its block matrix form with respect to the decomposition $\Lambda^k[\mathcal{X} \oplus \mathcal{Y}] \cong \Lambda^k[\mathcal{X}] \oplus \Lambda^k[\mathcal{Y}]$, we write
\begin{equation}\label{eq:deltaM_block}
    \delta_{\mathbf{M}}^k = \begin{pmatrix} \delta_{\mathbf{A}}^k & \delta_{\mathbf{C}}^k \\ 0 & \delta_{\mathbf{B}}^k \end{pmatrix}.
\end{equation}
Then, the identity $\delta_{\mathbf{M}}^{k+1} \delta_{\mathbf{M}}^k = 0$ is equivalent to the system
\begin{equation}\label{eq:system_expanded}
    \begin{cases}
        \delta_{\mathbf{A}}^{k+1} \delta_{\mathbf{A}}^k = 0 \\
        \delta_{\mathbf{B}}^{k+1} \delta_{\mathbf{B}}^k = 0 \\
        \delta_{\mathbf{A}}^{k+1} \delta_{\mathbf{C}}^k + \delta_{\mathbf{C}}^{k+1} \delta_{\mathbf{B}}^k = 0
    \end{cases}
\end{equation}
where the third equation is the lift of \eqref{eq:coupling_relation} to the exterior algebra $\Lambda[\mathcal{X} \oplus \mathcal{Y}]$.
\medskip

Taylor invertibility is defined by the purely algebraic exactness of the Koszul complex, but exactness forces every range $\mathcal{R}(\delta_{\mathbf{M}}^k) = \ker \delta_{\mathbf{M}}^{k+1}$ to be closed. The following theorem clarifies the extent to which this closedness is inherited by the diagonal boundary maps.
\begin{theorem}\label{thm:AB_range_closed}
    Let $\mathbf{M} = \begin{pmatrix} \mathbf{A} & \mathbf{C} \\ \mathbf{0} & \mathbf{B} \end{pmatrix}$ be the commuting $n$-tuple defined by \eqref{eq:M_tuple}. Then:
    \begin{enumerate}[\rm($i$)]
        \item If $\mathscr{K}(\mathbf{M})$ is exact at degree $0$, that is, $\ker \delta_{\mathbf{M}}^{0} = \{0\}$, then $\delta_{\mathbf{A}}^0$ is injective.
        \item If $\ker \delta_{\mathbf{M}}^{0} = \{0\}$ and $\mathcal{R}(\delta_{\mathbf{M}}^{0})$ is closed, which is the case when $\mathscr{K}(\mathbf{M})$ is exact at the degrees $0$ and $1$, then $\delta_{\mathbf{A}}^0$ is bounded below. In particular, $\mathcal{R}(\delta_{\mathbf{A}}^0)$ is closed in $\Lambda^1[\mathcal{X}]$.
        \item If $\mathscr{K}(\mathbf{M})$ is exact at degree $n$, that is, $\mathcal{R}(\delta_{\mathbf{M}}^{n-1}) = \Lambda^n[\mathcal{X} \oplus \mathcal{Y}]$, then $\delta_{\mathbf{B}}^{n-1}$ is surjective. In particular, $\mathcal{R}(\delta_{\mathbf{B}}^{n-1}) = \Lambda^n[\mathcal{Y}]$ is closed.
    \end{enumerate}
    Consequently, if $\mathbf{M}$ is Taylor invertible, then $\delta_{\mathbf{A}}^0$ is bounded below and $\delta_{\mathbf{B}}^{n-1}$ is surjective.
\end{theorem}
\begin{proof}
    ($i$) Let $x \in \ker \delta_{\mathbf{A}}^0$. Since
    \[
    \delta_{\mathbf{M}}^0 \begin{pmatrix} x \\ 0 \end{pmatrix} = \begin{pmatrix} \delta_{\mathbf{A}}^0 x \\ 0 \end{pmatrix} = 0,
    \]
    the vector $\begin{pmatrix} x \\ 0 \end{pmatrix}$ belongs to $\ker \delta_{\mathbf{M}}^0 = \{0\}$, whence $x = 0$.

    ($ii$) By hypothesis, $\delta_{\mathbf{M}}^0$ is injective and has closed range, so by the Bounded Inverse Theorem it is bounded below, meaning that there exists a constant $c > 0$ such that $\|\delta_{\mathbf{M}}^0 \zeta\| \ge c\|\zeta\|$ for all $\zeta \in \Lambda^0[\mathcal{X} \oplus \mathcal{Y}]$. Let $x \in \mathcal{X} \cong \Lambda^0[\mathcal{X}]$. Applying this to $\zeta = \begin{pmatrix} x \\ 0 \end{pmatrix}$ and using the identity displayed in ($i$), we obtain
    \[
    \|\delta_{\mathbf{A}}^0 x\| = \left\| \begin{pmatrix} \delta_{\mathbf{A}}^0 x \\ 0 \end{pmatrix} \right\| = \left\| \delta_{\mathbf{M}}^0 \begin{pmatrix} x \\ 0 \end{pmatrix} \right\| \ge c \left\| \begin{pmatrix} x \\ 0 \end{pmatrix} \right\| = c\|x\|.
    \]
    This implies that $\delta_{\mathbf{A}}^0$ is bounded below, and consequently, $\mathcal{R}(\delta_{\mathbf{A}}^0)$ is closed in $\Lambda^1[\mathcal{X}]$.

    ($iii$) Let $\eta \in \Lambda^n[\mathcal{Y}]$ be arbitrary. By hypothesis, the element $\begin{pmatrix} 0 \\ \eta \end{pmatrix} \in \Lambda^n[\mathcal{X} \oplus \mathcal{Y}]$ has a pre-image under $\delta_{\mathbf{M}}^{n-1}$. Let $\begin{pmatrix} \xi' \\ \eta' \end{pmatrix} \in \Lambda^{n-1}[\mathcal{X} \oplus \mathcal{Y}]$ be such that
    \[
    \delta_{\mathbf{M}}^{n-1} \begin{pmatrix} \xi' \\ \eta' \end{pmatrix} = \begin{pmatrix} \delta_{\mathbf{A}}^{n-1} \xi' + \delta_{\mathbf{C}}^{n-1} \eta' \\ \delta_{\mathbf{B}}^{n-1} \eta' \end{pmatrix} = \begin{pmatrix} 0 \\ \eta \end{pmatrix}.
    \]
    From the second component, we have $\delta_{\mathbf{B}}^{n-1} \eta' = \eta$. Since $\eta$ was arbitrary, it follows that $\delta_{\mathbf{B}}^{n-1}$ is surjective, which implies that $\mathcal{R}(\delta_{\mathbf{B}}^{n-1}) = \Lambda^n[\mathcal{Y}]$ is closed.

    Finally, if $\mathbf{M}$ is Taylor invertible, then $\mathscr{K}(\mathbf{M})$ is exact at every degree, so the hypotheses of ($ii$) and ($iii$) are satisfied.
\end{proof}
\begin{remark}\label{rem:intermediate_ranges}
    For $n = 1$, the invertibility of $M_C$ forces $A$ to be bounded below and $B$ to be surjective, so both diagonal entries automatically have closed ranges. Theorem \ref{thm:AB_range_closed} shows that, in the multivariable setting, the Taylor invertibility of $\mathbf{M}$ still controls the two \emph{extreme} positions: the bottom range of $\delta_{\mathbf{A}}$ and the top range of $\delta_{\mathbf{B}}$. The proof rests on the injectivity of $\delta_{\mathbf{M}}$ at the bottom, and its surjectivity at the top, of the complex $\mathscr{K}(\mathbf{M})$, and does not extend to the remaining degrees. This is not a defect of the method. Example \ref{ex:non_closed} below exhibits a Taylor invertible commuting pair of the form \eqref{eq:M_tuple} for which neither $\mathcal{R}(\delta_{\mathbf{B}}^{0})$ nor $\mathcal{R}(\delta_{\mathbf{A}}^{1})$ is closed. Let us stress that this failure is purely topological. As follows from the algebraic part of the long exact sequence (Lemma \ref{lem:zig_zag_Taylor}($i$) below), the Taylor invertibility of $\mathbf{M}$ always yields \emph{linear} isomorphisms $\partial_k \colon H^k(\mathbf{B}) \to H^{k+1}(\mathbf{A})$ in every degree, but this alone gives no control over the closedness of the individual ranges. This observation is upgraded to a full characterization in Theorem \ref{thm:taylor_char}.
\end{remark}

We now turn to the remaining degrees. In view of Remark \ref{rem:intermediate_ranges}, we look for sufficient conditions on the coupling tuple $\mathbf{C}$ under which the exactness, or the closedness of the ranges, of $\delta_{\mathbf{M}}$ transfers to the diagonal boundary maps. The results that follow complement Theorem \ref{thm:AB_range_closed}.

\begin{theorem}\label{thm:B_range_closed}
    Let $\mathbf{M}$ be the commuting $n$-tuple defined by \eqref{eq:M_tuple}. Fix $0 \leq k < n-1$ and assume that $\mathcal{R}(\delta_{\mathbf{M}}^k) = \ker(\delta_{\mathbf{M}}^{k+1})$. If there exists a constant $L \geq 0$ such that
    \begin{equation}\label{eq:domination}
        \|\delta_{\mathbf{C}}^k\eta\| \leq L\,\|\delta_{\mathbf{B}}^k\eta\| \qquad \text{for all } \eta \in \Lambda^k[\mathcal{Y}],
    \end{equation}
    then $\mathcal{R}(\delta_{\mathbf{B}}^k)$ is closed.
\end{theorem}
Let us note that \eqref{eq:domination} holds if and only if there exists a bounded operator $T$ from the closure of $\mathcal{R}(\delta_{\mathbf{B}}^k)$ into $\Lambda^{k+1}[\mathcal{X}]$ such that $\delta_{\mathbf{C}}^k = T \delta_{\mathbf{B}}^k$. Indeed, under \eqref{eq:domination} the assignment $\delta_{\mathbf{B}}^k \eta \mapsto \delta_{\mathbf{C}}^k \eta$ is well defined and bounded on $\mathcal{R}(\delta_{\mathbf{B}}^k)$ and extends by continuity to its closure, while the converse implication is obvious.
\begin{proof}
    Let $\eta$ be in the closure of $\mathcal{R}(\delta_{\mathbf{B}}^k)$. There exists a sequence $(\eta_m)_{m \geq 1}$ in $\Lambda^k[\mathcal{Y}]$ with $\delta_{\mathbf{B}}^k\eta_m \to \eta$. In particular, $(\delta_{\mathbf{B}}^k\eta_m)_{m \geq 1}$ is a Cauchy sequence, and by \eqref{eq:domination},
    \[
    \|\delta_{\mathbf{C}}^k\eta_m - \delta_{\mathbf{C}}^k\eta_l\| = \|\delta_{\mathbf{C}}^k(\eta_m - \eta_l)\| \leq L\,\|\delta_{\mathbf{B}}^k(\eta_m - \eta_l)\| \longrightarrow 0 \quad (m, l \to \infty).
    \]
    Hence $(\delta_{\mathbf{C}}^k\eta_m)_{m \geq 1}$ is a Cauchy sequence in the Banach space $\Lambda^{k+1}[\mathcal{X}]$, and therefore converges to some $\xi \in \Lambda^{k+1}[\mathcal{X}]$.

    We compute, using the third relation in \eqref{eq:system_expanded}:
    \[
    \delta_{\mathbf{A}}^{k+1}\xi + \delta_{\mathbf{C}}^{k+1}\eta = \lim_{m \to \infty} \Bigl(\delta_{\mathbf{A}}^{k+1}\delta_{\mathbf{C}}^k\eta_{m} + \delta_{\mathbf{C}}^{k+1}\delta_{\mathbf{B}}^k\eta_{m}\Bigr) = \lim_{m \to \infty} 0 = 0,
    \]
    and $\delta_{\mathbf{B}}^{k+1}\eta = \lim_{m \to \infty} \delta_{\mathbf{B}}^{k+1}\delta_{\mathbf{B}}^k\eta_{m} = 0$. Therefore,
    \[
    \delta_{\mathbf{M}}^{k+1}\begin{pmatrix} \xi \\ \eta \end{pmatrix} = \begin{pmatrix} \delta_{\mathbf{A}}^{k+1}\xi + \delta_{\mathbf{C}}^{k+1}\eta \\ \delta_{\mathbf{B}}^{k+1}\eta \end{pmatrix} = \begin{pmatrix} 0 \\ 0 \end{pmatrix},
    \]
    so $\begin{pmatrix} \xi \\ \eta \end{pmatrix} \in \ker(\delta_{\mathbf{M}}^{k+1}) = \mathcal{R}(\delta_{\mathbf{M}}^k)$. Hence there exists $\begin{pmatrix} \tilde\xi \\ \tilde\eta \end{pmatrix}$ with $\delta_{\mathbf{M}}^k \begin{pmatrix} \tilde\xi \\ \tilde\eta \end{pmatrix} = \begin{pmatrix} \xi \\ \eta \end{pmatrix}$. Reading off the second component, $\eta = \delta_{\mathbf{B}}^k\tilde\eta \in \mathcal{R}(\delta_{\mathbf{B}}^k)$.
\end{proof}

\begin{corollary}\label{cor:B_range_h}
    Let $\mathbf{M}$ be the commuting $n$-tuple defined by \eqref{eq:M_tuple}. Fix $0 \leq k < n-1$ and assume that $\mathcal{R}(\delta_{\mathbf{M}}^k) = \ker(\delta_{\mathbf{M}}^{k+1})$. Then $\mathcal{R}(\delta_{\mathbf{B}}^k)$ is closed if either of the following conditions holds:
    \begin{enumerate}[\rm($i$)]
        \item There exists $D \in \mathfrak{B}(\mathcal{Y}, \mathcal{X})$ such that $C_i = DB_i$ for all $1 \leq i \leq n$.
        \item $k = 0$ and there exists a constant $M \geq 0$ such that $\|C_i y\| \leq M \max_{1 \leq j \leq n}\|B_j y\|$ for all $y \in \mathcal{Y}$ and all $i$.
    \end{enumerate}
\end{corollary}

\begin{proof}
    ($i$) Extending $D$ diagonally to the exterior algebra as in Section \ref{sec:prelim}, we have $\delta_{\mathbf{C}}^k = D \circ \delta_{\mathbf{B}}^k$; indeed, for $\eta = \sum_{|F|=k} y_F \otimes e_F$,
    \[
    \delta_{\mathbf{C}}^k \eta = \sum_{F} \sum_{i=1}^n C_i y_F \otimes (e_i \wedge e_F) = \sum_{F} \sum_{i=1}^n D B_i y_F \otimes (e_i \wedge e_F) = D\bigl(\delta_{\mathbf{B}}^k \eta\bigr).
    \]
    Hence $\|\delta_{\mathbf{C}}^k\eta\| \leq \|D\|\,\|\delta_{\mathbf{B}}^k\eta\|$, which is \eqref{eq:domination} with $L = \|D\|$, and Theorem \ref{thm:B_range_closed} applies.

    ($ii$) For $k = 0$ and $\eta = y \in \Lambda^0[\mathcal{Y}] \cong \mathcal{Y}$, the norm conventions of Section \ref{sec:prelim} give $\|\delta_{\mathbf{B}}^0 y\| = \max_{j}\|B_j y\|$ and $\|\delta_{\mathbf{C}}^0 y\| = \max_{i}\|C_i y\| \leq M \max_{j}\|B_j y\| = M \|\delta_{\mathbf{B}}^0 y\|$, which is \eqref{eq:domination} with $L = M$, and Theorem \ref{thm:B_range_closed} applies again.
\end{proof}

\begin{theorem}\label{thm:A_range_closed}
    Let $\mathbf{M}$ be the commuting $n$-tuple defined by \eqref{eq:M_tuple}. Fix $1 \leq k \leq n-1$ and assume that $\mathcal{R}(\delta_{\mathbf{M}}^k)$ is closed and that $\ker(\delta_{\mathbf{B}}^k) \subseteq \ker(\delta_{\mathbf{C}}^k)$. Then $\mathcal{R}(\delta_{\mathbf{A}}^k)$ is closed.
\end{theorem}

\begin{proof}
    Let $\xi$ be in the closure of $\mathcal{R}(\delta_{\mathbf{A}}^k)$. Then $\begin{pmatrix} \xi \\ 0 \end{pmatrix}$ lies in the closure of $\mathcal{R}(\delta_{\mathbf{M}}^k)$. Since $\mathcal{R}(\delta_{\mathbf{M}}^k)$ is closed by hypothesis, there exists $\begin{pmatrix} \tilde\xi \\ \tilde\eta \end{pmatrix}$ such that
    \[
    \delta_{\mathbf{M}}^k\begin{pmatrix} \tilde\xi \\ \tilde\eta \end{pmatrix} = \begin{pmatrix} \delta_{\mathbf{A}}^k\tilde\xi + \delta_{\mathbf{C}}^k\tilde\eta \\ \delta_{\mathbf{B}}^k\tilde\eta \end{pmatrix} = \begin{pmatrix} \xi \\ 0 \end{pmatrix}.
    \]
    The second component gives $\delta_{\mathbf{B}}^k\tilde\eta = 0$, i.e., $\tilde\eta \in \ker(\delta_{\mathbf{B}}^k)$. By hypothesis, $\ker(\delta_{\mathbf{B}}^k) \subseteq \ker(\delta_{\mathbf{C}}^k)$, so $\delta_{\mathbf{C}}^k\tilde\eta = 0$. The first component then gives $\xi = \delta_{\mathbf{A}}^k\tilde\xi \in \mathcal{R}(\delta_{\mathbf{A}}^k)$.
\end{proof}

\begin{corollary}\label{cor:A_range_h}
    Let $\mathbf{M}$ be the commuting $n$-tuple defined by \eqref{eq:M_tuple}, fix $1 \leq k \leq n-1$, and assume that $\mathcal{R}(\delta_{\mathbf{M}}^k)$ is closed. If \eqref{eq:domination} holds for some constant $L \geq 0$, in particular if $C_i = D B_i$ for all $1 \leq i \leq n$ and some $D \in \mathfrak{B}(\mathcal{Y}, \mathcal{X})$, then $\mathcal{R}(\delta_{\mathbf{A}}^k)$ is closed.
\end{corollary}

\begin{proof}
    By Theorem \ref{thm:A_range_closed}, it suffices to verify the inclusion $\ker(\delta_{\mathbf{B}}^k) \subseteq \ker(\delta_{\mathbf{C}}^k)$, which is immediate from \eqref{eq:domination}. If $C_i = D B_i$ for all $i$, then \eqref{eq:domination} holds with $L = \|D\|$, as shown in the proof of Corollary \ref{cor:B_range_h}.
\end{proof}

\begin{remark}\label{rem:componentwise}
    The restriction to $k = 0$ in Corollary \ref{cor:B_range_h}($ii$) cannot be removed: for $k \geq 1$, componentwise domination of $\mathbf{C}$ by $\mathbf{B}$ no longer controls $\delta_{\mathbf{C}}^k$ by $\delta_{\mathbf{B}}^k$, due to cancellations inside the exterior powers. For instance, let $n = 2$, $\mathcal{X} = \mathcal{Y}$, let $B \in \mathfrak{B}(\mathcal{Y})$ be a nonzero operator with $B^2 = 0$, and put $B_1 = B_2 = B$, $C_1 = B$, $C_2 = -B$, $A_1 = A_2 = 0$. Then \eqref{eq:coupling_relation} holds, $\mathbf{M}$ is a commuting pair, and $\|C_i y\| \leq \max_{j}\|B_j y\|$ for all $y \in \mathcal{Y}$, while for $\eta = y \otimes e_1 + y \otimes e_2$ we have
    \[
    \delta_{\mathbf{B}}^1 \eta = 0 \qquad \text{and} \qquad \delta_{\mathbf{C}}^1 \eta = 2By \otimes (e_1 \wedge e_2) \neq 0
    \]
    whenever $y \notin \mathcal{N}(B)$. In particular, $\ker \delta_{\mathbf{B}}^1 \not\subseteq \ker \delta_{\mathbf{C}}^1$, so that \eqref{eq:domination} fails at degree $k = 1$ for every $L \geq 0$.
\end{remark}

The following example shows that the sufficient conditions obtained above cannot be dispensed with: the Taylor invertibility of $\mathbf{M}$ alone does not force the closedness of the diagonal ranges in the intermediate degrees.

\begin{example}\label{ex:non_closed}
    Let $\mathcal{H}$ be an infinite-dimensional Hilbert space and let $K \in \mathfrak{B}(\mathcal{H})$ be an operator whose range is not closed, for instance $K = \operatorname{diag}(1, \tfrac{1}{2}, \tfrac{1}{3}, \dots)$ with respect to an orthonormal basis of a separable $\mathcal{H}$. Let $\mathcal{Z} = \ell^2(\mathbb{Z}; \mathcal{H})$ be the Hilbert space of all square-summable sequences $z = (z_i)_{i \in \mathbb{Z}}$ with $z_i \in \mathcal{H}$, and define $M_1, M_2 \in \mathfrak{B}(\mathcal{Z})$ by
    \[
    (M_1 z)_i = z_{i+1} \qquad \text{and} \qquad (M_2 z)_i = K z_i \qquad (i \in \mathbb{Z}).
    \]
    Then $M_1$ is unitary, $M_1$ and $M_2$ commute, and the closed subspace
    \[
    \mathcal{X} = \{ z \in \mathcal{Z} : z_i = 0 \text{ for all } i \geq 0 \}
    \]
    is invariant under both $M_1$ and $M_2$. Put $\mathcal{Y} = \mathcal{X}^{\perp} = \{ z \in \mathcal{Z} : z_i = 0 \text{ for all } i < 0 \}$, which we identify with $\ell^2(\mathbb{N}_0; \mathcal{H})$. With respect to $\mathcal{Z} = \mathcal{X} \oplus \mathcal{Y}$, the commuting pair $\mathbf{M} = (M_1, M_2)$ is of the form \eqref{eq:M_tuple}, and its compressions to $\mathcal{Y}$ are given by
    \[
    (B_1 y)_i = y_{i+1} \qquad \text{and} \qquad (B_2 y)_i = K y_i \qquad (i \in \mathbb{N}_0).
    \]
    The Hilbert space norm on $\mathcal{Z}$ is equivalent to the norm of Section \ref{sec:prelim}, and none of the notions involved is affected by passing to an equivalent norm. Since $M_1$ is invertible, $\mathbf{M}$ is Taylor invertible, because a commuting tuple one of whose entries is invertible is Taylor invertible (see \cite{Taylor70}).

    We claim that $\mathcal{R}(\delta_{\mathbf{B}}^0)$ is not closed. Let $V \colon \mathcal{H} \to \mathcal{Y}$ be the isometry $Vh = (h, 0, 0, \dots)$. Then $B_1 V h = 0$ and $B_2 V h = V K h$ for every $h \in \mathcal{H}$, so that
    \[
    \delta_{\mathbf{B}}^0 (Vh) = VKh \otimes e_2 \qquad \text{and} \qquad \|\delta_{\mathbf{B}}^0 (Vh)\| = \|Kh\|.
    \]
    Moreover, $\ker B_1 = V(\mathcal{H})$ and $Vh \in \ker B_2$ if and only if $h \in \ker K$, whence $\ker \delta_{\mathbf{B}}^0 = \ker B_1 \cap \ker B_2 = V(\ker K)$ and $\operatorname{dist}(Vh, \ker \delta_{\mathbf{B}}^0) = \operatorname{dist}(h, \ker K)$. If $\mathcal{R}(\delta_{\mathbf{B}}^0)$ were closed, then the operator induced by $\delta_{\mathbf{B}}^0$ on $\Lambda^0[\mathcal{Y}]/\ker \delta_{\mathbf{B}}^0$ would be a bounded bijection onto the Banach space $\mathcal{R}(\delta_{\mathbf{B}}^0)$, and the Bounded Inverse Theorem would provide a constant $c > 0$ such that
    \[
    \|\delta_{\mathbf{B}}^0 \eta\| \geq c \operatorname{dist}(\eta, \ker \delta_{\mathbf{B}}^0) \qquad \text{for all } \eta \in \Lambda^0[\mathcal{Y}].
    \]
    Applied to $\eta = Vh$, this gives $\|Kh\| \geq c \operatorname{dist}(h, \ker K)$ for all $h \in \mathcal{H}$, which forces $\mathcal{R}(K)$ to be closed, a contradiction. Hence $\mathcal{R}(\delta_{\mathbf{B}}^0)$ is not closed, although $\mathbf{M}$ is Taylor invertible.

    The range $\mathcal{R}(\delta_{\mathbf{A}}^1)$ is not closed either. Identifying $\mathcal{X}$ with $\ell^2(\mathbb{N}_0; \mathcal{H})$ via $x_j = z_{-1-j}$ ($j \in \mathbb{N}_0$), the restrictions of $M_1$ and $M_2$ to $\mathcal{X}$ become
    \[
    (A_1 x)_0 = 0, \qquad (A_1 x)_j = x_{j-1} \quad (j \geq 1), \qquad (A_2 x)_j = K x_j \quad (j \in \mathbb{N}_0).
    \]
    Under the identification $\Lambda^2[\mathcal{X}] \cong \mathcal{X}$, the boundary map $\delta_{\mathbf{A}}^1$ acts as $(x^{(1)}, x^{(2)}) \mapsto A_1 x^{(2)} - A_2 x^{(1)}$, so that $\mathcal{R}(\delta_{\mathbf{A}}^1) = \mathcal{R}(A_1) + \mathcal{R}(A_2)$. Since $\mathcal{R}(A_1)$ consists of all sequences with vanishing first term, we obtain
    \[
    \mathcal{R}(\delta_{\mathbf{A}}^1) = \{ x \in \ell^2(\mathbb{N}_0; \mathcal{H}) : x_0 \in \mathcal{R}(K) \},
    \]
    whose closure is $\{ x \in \ell^2(\mathbb{N}_0; \mathcal{H}) : x_0 \in \overline{\mathcal{R}(K)} \}$. As $\mathcal{R}(K)$ is not closed, neither is $\mathcal{R}(\delta_{\mathbf{A}}^1)$.
\end{example}

\begin{remark}\label{rem:non_closed_n}
    Example \ref{ex:non_closed} shows that already for $n = 2$ both $\mathcal{R}(\delta_{\mathbf{B}}^0)$ and $\mathcal{R}(\delta_{\mathbf{A}}^1)$ may fail to be closed for a Taylor invertible $\mathbf{M}$. Adjoining zero operators, that is, passing from $(M_1, M_2)$ to $(M_1, M_2, 0, \dots, 0)$, produces Taylor invertible commuting $n$-tuples with the same defect in every intermediate degree. Indeed, the Koszul complex of the enlarged tuple is the tensor product of the Koszul complex of $(M_1, M_2)$ with the exterior algebra on the additional $n-2$ generators, equipped with the zero differential, so that $\mathcal{R}(\delta_{\mathbf{B}}^k)$ is a finite direct sum of copies of $\mathcal{R}(\delta_{(B_1, B_2)}^0)$ and $\mathcal{R}(\delta_{(B_1, B_2)}^1)$, and similarly for $\mathbf{A}$. Consequently, $\mathcal{R}(\delta_{\mathbf{B}}^k)$ fails to be closed for every $0 \leq k \leq n-2$, and $\mathcal{R}(\delta_{\mathbf{A}}^k)$ fails to be closed for every $1 \leq k \leq n-1$. In this sense Theorem \ref{thm:AB_range_closed} is sharp, and the closedness of the remaining ranges has to be imposed as a hypothesis or ensured by conditions such as those of Theorems \ref{thm:B_range_closed} and \ref{thm:A_range_closed}. Let us also note that in \cite{Curto80} the diagonal tuples are a priori Taylor invertible, so that all their boundary maps have closed ranges and this phenomenon does not arise there.
\end{remark}

\begin{theorem}\label{thm:M_range_closed}
    Let $\mathbf{M}$ be the commuting $n$-tuple defined by \eqref{eq:M_tuple}. Fix $k \in \mathbb{Z}$ and assume that the boundary maps $\delta_{\mathbf{A}}^k$ and $\delta_{\mathbf{B}}^k$ have closed ranges. Then the range of $\delta_{\mathbf{M}}^k$ is closed in $\Lambda^{k+1}[\mathcal{X} \oplus \mathcal{Y}]$ if and only if the induced linear map
    \[
    \partial_k \colon H^k(\mathbf{B}) \to H^{k+1}(\mathbf{A}), \qquad \partial_k([\eta]) = [\delta_{\mathbf{C}}^k \eta],
    \]
    has closed range in $H^{k+1}(\mathbf{A})$.
\end{theorem}
Let us note that, as will be seen in the proof, the closedness of $\mathcal{R}(\partial_k)$ is equivalent to the closedness of the subspace $\delta_{\mathbf{C}}^k\bigl(\ker \delta_{\mathbf{B}}^k\bigr) + \mathcal{R}(\delta_{\mathbf{A}}^k)$ in $\Lambda^{k+1}[\mathcal{X}]$.
\begin{proof}
    Before proceeding to the proof of closedness, we first must ensure that the induced map $\partial_k$ is well-defined. By the commutativity of $\mathbf{M}$, specifically the third relation in \eqref{eq:system_expanded}, we have that for any cocycle $\eta \in \ker \delta_{\mathbf{B}}^k$,
    \begin{equation*}
        \delta_{\mathbf{A}}^{k+1} (\delta_{\mathbf{C}}^k \eta) = - \delta_{\mathbf{C}}^{k+1} (\delta_{\mathbf{B}}^k \eta) = 0,
    \end{equation*}
    which confirms that $\delta_{\mathbf{C}}^k \eta \in \ker \delta_{\mathbf{A}}^{k+1}$. Thus, the assignment $\partial_k([\eta]) = [\delta_{\mathbf{C}}^k \eta]$ correctly defines a mapping between the cohomology groups $H^k(\mathbf{B})$ and $H^{k+1}(\mathbf{A})$.
    Next, let us show that the definition is independent of the choice of the representative $\eta \in [\eta]$. Let $\eta_1, \eta_2 \in \ker \delta_{\mathbf{B}}^k$ be two representatives of the same cohomology class, such that $[\eta_1] = [\eta_2]$. This implies that $\eta_1 - \eta_2 \in \mathcal{R}(\delta_{\mathbf{B}}^{k-1})$, meaning there exists some $\zeta \in \Lambda^{k-1}[\mathcal{Y}]$ such that $\eta_1 - \eta_2 = \delta_{\mathbf{B}}^{k-1} \zeta$.
    Applying the coupling operator $\delta_{\mathbf{C}}^k$ to this difference, we have
    \begin{equation*}
        \delta_{\mathbf{C}}^k (\eta_1 - \eta_2) = \delta_{\mathbf{C}}^k (\delta_{\mathbf{B}}^{k-1} \zeta)=- \delta_{\mathbf{A}}^k (\delta_{\mathbf{C}}^{k-1} \zeta).
    \end{equation*}
    Since $\delta_{\mathbf{A}}^k (\delta_{\mathbf{C}}^{k-1} \zeta) \in \mathcal{R}(\delta_{\mathbf{A}}^k)$, it follows that $\delta_{\mathbf{C}}^k \eta_1$ and $\delta_{\mathbf{C}}^k \eta_2$ differ only by an element of $\mathcal{R}(\delta_{\mathbf{A}}^k)$, and thus
    \begin{equation*}
        [\delta_{\mathbf{C}}^k \eta_1] = [\delta_{\mathbf{C}}^k \eta_2] \in H^{k+1}(\mathbf{A}).
    \end{equation*}
    This confirms that the map $\partial_k$ is well-defined.
    \medskip

    Next, we fix the topological framework. Since $\mathcal{R}(\delta_{\mathbf{A}}^k)$ is closed, the quotient
    \[
    Q_{k+1} := \Lambda^{k+1}[\mathcal{X}] / \mathcal{R}(\delta_{\mathbf{A}}^k)
    \]
    is a Banach space; let $\pi \colon \Lambda^{k+1}[\mathcal{X}] \to Q_{k+1}$ denote the canonical (bounded) quotient projection. The subspace $\ker \delta_{\mathbf{A}}^{k+1}$ is closed and contains $\mathcal{R}(\delta_{\mathbf{A}}^k)$, so $H^{k+1}(\mathbf{A}) = \pi(\ker \delta_{\mathbf{A}}^{k+1})$ is a closed subspace of $Q_{k+1}$, and $\mathcal{R}(\partial_k) \subseteq H^{k+1}(\mathbf{A})$. Consequently, $\mathcal{R}(\partial_k)$ is closed in $H^{k+1}(\mathbf{A})$ if and only if it is closed in $Q_{k+1}$, which, in turn, is equivalent to the closedness of $\pi^{-1}(\mathcal{R}(\partial_k)) = \delta_{\mathbf{C}}^k(\ker \delta_{\mathbf{B}}^k) + \mathcal{R}(\delta_{\mathbf{A}}^k)$ in $\Lambda^{k+1}[\mathcal{X}]$. This justifies the equivalent formulation stated before the proof.
    \medskip

    Assume first that $\mathcal{R}(\partial_k)$ is closed, and let $\begin{pmatrix} u \\ v \end{pmatrix}$ belong to the closure of $\mathcal{R}(\delta_{\mathbf{M}}^k)$. Then there exists a sequence $\begin{pmatrix} x_m \\ y_m \end{pmatrix}_{m \geq 1}$ in $\Lambda^k[\mathcal{X}] \oplus \Lambda^k[\mathcal{Y}]$ such that
    \begin{equation}\label{eq:m_conv_1}
        \delta_{\mathbf{A}}^k x_m + \delta_{\mathbf{C}}^k y_m \to u
    \end{equation}
    and
    \begin{equation}\label{eq:m_conv_2}
        \delta_{\mathbf{B}}^k y_m \to v.
    \end{equation}
    Since $\mathcal{R}(\delta_{\mathbf{B}}^k)$ is closed by hypothesis, \eqref{eq:m_conv_2} implies that $v \in \mathcal{R}(\delta_{\mathbf{B}}^k)$, hence there exists some $y_0 \in \Lambda^k[\mathcal{Y}]$ such that $\delta_{\mathbf{B}}^k y_0 = v$.
    Now, consider the operator $\delta_{\mathbf{B}}^k$ acting from $\Lambda^k[\mathcal{Y}]$ to $\mathcal{R}(\delta_{\mathbf{B}}^k)$. By the Bounded Inverse Theorem, the induced operator $\widehat{\delta_{\mathbf{B}}^k}: \Lambda^k[\mathcal{Y}]/\ker \delta_{\mathbf{B}}^k \to \mathcal{R}(\delta_{\mathbf{B}}^k)$ is a topological isomorphism. Since $\delta_{\mathbf{B}}^k(y_m - y_0) \to 0$, the sequence of cosets $[y_m - y_0]$ converges to zero in the quotient space $\Lambda^k[\mathcal{Y}]/\ker \delta_{\mathbf{B}}^k$.
    By the definition of the quotient norm, for every $m$ we can choose $\theta_m \in \ker \delta_{\mathbf{B}}^k$ with $\|y_m - y_0 - \theta_m\| \leq \|[y_m - y_0]\| + 1/m$. Putting $z_m = y_m - y_0 - \theta_m$, we obtain the decomposition $y_m - y_0 = \theta_m + z_m$, where $\theta_m \in \ker \delta_{\mathbf{B}}^k$ and $z_m \to 0$ in $\Lambda^k[\mathcal{Y}]$. Substituting $y_m = y_0 + \theta_m + z_m$ into \eqref{eq:m_conv_1}, we obtain
    \[ \delta_{\mathbf{A}}^k x_m + \delta_{\mathbf{C}}^k (y_0 + \theta_m + z_m) \to u. \]
    Since $\delta_{\mathbf{C}}^k$ is bounded and $z_m \to 0$, the term $\delta_{\mathbf{C}}^k z_m$ vanishes in the limit. Defining $w = u - \delta_{\mathbf{C}}^k y_0$, we have
    \begin{equation}\label{eq:aux_w}
        \delta_{\mathbf{A}}^k x_m + \delta_{\mathbf{C}}^k \theta_m \to w.
    \end{equation}
    Applying the (continuous) projection $\pi$ to \eqref{eq:aux_w} and noting that $\pi(\delta_{\mathbf{A}}^k x_m) = 0$, we get
    \[ \pi(\delta_{\mathbf{C}}^k \theta_m) = \partial_k([\theta_m]) \to \pi(w) \quad \text{in } Q_{k+1}. \]
    By our assumption and the discussion above, $\mathcal{R}(\partial_k)$ is closed in $Q_{k+1}$, hence the limit $\pi(w)$ belongs to $\mathcal{R}(\partial_k)$. This implies there exists some $[\kappa] \in H^k(\mathbf{B})$ (represented by $\kappa \in \ker \delta_{\mathbf{B}}^k$) such that $\partial_k([\kappa]) = \pi(w)$. This is equivalent to stating that $w - \delta_{\mathbf{C}}^k \kappa \in \mathcal{R}(\delta_{\mathbf{A}}^k)$.
    Thus, there exists $x \in \Lambda^k[\mathcal{X}]$ such that $\delta_{\mathbf{A}}^k x = w - \delta_{\mathbf{C}}^k \kappa$. Finally, by choosing $y = y_0 + \kappa \in \Lambda^k[\mathcal{Y}]$, we verify that:
    \[ \delta_{\mathbf{A}}^k x + \delta_{\mathbf{C}}^k y = (w - \delta_{\mathbf{C}}^k \kappa) + \delta_{\mathbf{C}}^k y_0 + \delta_{\mathbf{C}}^k \kappa = w + \delta_{\mathbf{C}}^k y_0 = u, \]
    and
    \[ \delta_{\mathbf{B}}^k y = \delta_{\mathbf{B}}^k (y_0 + \kappa) = \delta_{\mathbf{B}}^k y_0 = v. \]
    Therefore, $\delta_{\mathbf{M}}^k \begin{pmatrix} x \\ y \end{pmatrix} = \begin{pmatrix} u \\ v \end{pmatrix}$, which proves that $\mathcal{R}(\delta_{\mathbf{M}}^k)$ is closed.
    \medskip

    Conversely, assume that $\mathcal{R}(\delta_{\mathbf{M}}^k)$ is closed. An element
    \[
    \delta_{\mathbf{M}}^k \begin{pmatrix} \xi \\ \eta \end{pmatrix} = \begin{pmatrix} \delta_{\mathbf{A}}^k \xi + \delta_{\mathbf{C}}^k \eta \\ \delta_{\mathbf{B}}^k \eta \end{pmatrix}
    \]
    of $\mathcal{R}(\delta_{\mathbf{M}}^k)$ belongs to the closed subspace $\Lambda^{k+1}[\mathcal{X}] \oplus \{0\}$ if and only if $\eta \in \ker \delta_{\mathbf{B}}^k$. Hence
    \[
    \mathcal{R}(\delta_{\mathbf{M}}^k) \cap \bigl(\Lambda^{k+1}[\mathcal{X}] \oplus \{0\}\bigr) = \bigl(\mathcal{R}(\delta_{\mathbf{A}}^k) + \delta_{\mathbf{C}}^k(\ker \delta_{\mathbf{B}}^k)\bigr) \oplus \{0\}.
    \]
    The left-hand side is closed as the intersection of two closed subspaces. Since the map $\xi \mapsto \begin{pmatrix} \xi \\ 0 \end{pmatrix}$ is an isometry, the subspace $\mathcal{R}(\delta_{\mathbf{A}}^k) + \delta_{\mathbf{C}}^k(\ker \delta_{\mathbf{B}}^k)$ is closed in $\Lambda^{k+1}[\mathcal{X}]$, and by the discussion above $\mathcal{R}(\partial_k)$ is closed.
\end{proof}
    \medskip
  We define the natural inclusion $i_k: \Lambda^k[\mathcal{X}] \to \Lambda^k[\mathcal{X} \oplus \mathcal{Y}]$ and the canonical projection $p_k: \Lambda^k[\mathcal{X} \oplus \mathcal{Y}] \to \Lambda^k[\mathcal{Y}]$ by
    \begin{equation*}
        i_k(\xi) = \begin{pmatrix} \xi \\ 0 \end{pmatrix} \quad \text{and} \quad p_k \begin{pmatrix} \xi \\ \eta \end{pmatrix} = \eta.
    \end{equation*}
    These maps are bounded linear operators and satisfy $\mathcal{R}(i_k) = \mathcal{N}(p_k)$, which ensures that the rows of the following diagram are exact:
    \begin{equation}\label{eq:koszul_chains}
        \begin{tikzcd}[row sep=1.2cm, column sep=1.3cm]
            & \vdots \arrow[d] & \vdots \arrow[d] & \vdots \arrow[d] & \\
            0 \arrow[r] & \Lambda^{k-1}[\mathcal{X}] \arrow[r, "i_{k-1}"] \arrow[d, "\delta_{\mathbf{A}}^{k-1}"] & \Lambda^{k-1}[\mathcal{X} \oplus \mathcal{Y}] \arrow[r, "p_{k-1}"] \arrow[d, "\delta_{\mathbf{M}}^{k-1}"] & \Lambda^{k-1}[\mathcal{Y}] \arrow[r] \arrow[d, "\delta_{\mathbf{B}}^{k-1}"] & 0 \\
            0 \arrow[r] & \Lambda^k[\mathcal{X}] \arrow[r, "i_k"] \arrow[d, "\delta_{\mathbf{A}}^k"] & \Lambda^k[\mathcal{X} \oplus \mathcal{Y}] \arrow[r, "p_k"] \arrow[d, "\delta_{\mathbf{M}}^k"] & \Lambda^k[\mathcal{Y}] \arrow[r] \arrow[d, "\delta_{\mathbf{B}}^k"] & 0 \\
            0 \arrow[r] & \Lambda^{k+1}[\mathcal{X}] \arrow[r, "i_{k+1}"] \arrow[d, "\delta_{\mathbf{A}}^{k+1}"] & \Lambda^{k+1}[\mathcal{X} \oplus \mathcal{Y}] \arrow[r, "p_{k+1}"] \arrow[d, "\delta_{\mathbf{M}}^{k+1}"] & \Lambda^{k+1}[\mathcal{Y}] \arrow[r] \arrow[d, "\delta_{\mathbf{B}}^{k+1}"] & 0 \\
            & \vdots & \vdots & \vdots &
        \end{tikzcd}
    \end{equation}
    The following lemma establishes that the boundary operators are compatible with the inclusion and projection maps, ensuring the commutativity of the diagram.
\begin{lemma}\label{lem:comm_diag}
    The Koszul boundary operators $\delta_{\mathbf{A}}, \delta_{\mathbf{B}}$, and $\delta_{\mathbf{M}}$ satisfy the following commutativity relations for every $k \in \{0, \dots, n-1\}$:
    \begin{enumerate}[\rm($i$)]
        \item $\delta_{\mathbf{M}}^k \circ i_k = i_{k+1} \circ \delta_{\mathbf{A}}^k$,
        \item $p_{k+1} \circ \delta_{\mathbf{M}}^k = \delta_{\mathbf{B}}^k \circ p_k$.
    \end{enumerate}
\end{lemma}
\begin{proof}
    Both identities follow immediately from the block representation \eqref{eq:deltaM_block}.
\end{proof}
    \begin{lemma}\label{lem:zig_zag_Taylor}
    Let $\mathbf{M} = \begin{pmatrix} \mathbf{A} & \mathbf{C} \\ \mathbf{0} & \mathbf{B} \end{pmatrix}$ be the commuting $n$-tuple of operators on $\mathcal{X} \oplus \mathcal{Y}$ defined in \eqref{eq:M_tuple}.
    \begin{enumerate}[\rm($i$)]
        \item The short exact sequence of Koszul complexes \eqref{eq:koszul_chains} induces a long exact sequence of vector spaces and linear maps
    \begin{equation}\label{eq:zig_zag_Taylor}
        \begin{tikzcd}[row sep=3.3em, column sep=large]
            & \vdots \arrow[dl] & \\
            H^k(\mathbf{A}) \arrow[r, "\iota_k^*"] & H^k(\mathbf{M}) \arrow[r, "\pi_k^*"] & H^k(\mathbf{B}) \arrow[dll, "\partial_k"'] \\
            H^{k+1}(\mathbf{A}) \arrow[r, "\iota_{k+1}^*"] & H^{k+1}(\mathbf{M}) \arrow[r, "\pi_{k+1}^*"] & H^{k+1}(\mathbf{B}) \arrow[dll, "\partial_{k+1}"'] \\
            H^{k+2}(\mathbf{A}) \arrow[r, "\iota_{k+2}^*"] & H^{k+2}(\mathbf{M}) \arrow[r, "\pi_{k+2}^*"] & H^{k+2}(\mathbf{B}) \arrow[dl] \\
            & \vdots &
        \end{tikzcd}
    \end{equation}
    where $\iota_k^*$ and $\pi_k^*$ denote the natural maps induced by the canonical injection $i_k$ and projection $p_k$, respectively, while the connecting homomorphism $H^k(\mathbf{B}) \to H^{k+1}(\mathbf{A})$ is given by $[\eta] \mapsto [\delta_{\mathbf{C}}^k \eta]$, that is, it is precisely the map $\partial_k$ of Theorem \ref{thm:M_range_closed}. No closedness assumptions are needed for this part.
        \item If, in addition, the boundary maps $\delta_{\mathbf{A}}^k$, $\delta_{\mathbf{B}}^k$, and $\delta_{\mathbf{M}}^k$ have closed ranges for all $k \in \mathbb{Z}$, then all the cohomology spaces in \eqref{eq:zig_zag_Taylor} are Banach spaces and all the maps in \eqref{eq:zig_zag_Taylor} are bounded; in particular,
        \[
        \|\partial_k [\eta]\| \leq \|\delta_{\mathbf{C}}^k\|\,\|[\eta]\| \qquad \text{for every } [\eta] \in H^k(\mathbf{B}).
        \]
    \end{enumerate}
\end{lemma}
\begin{proof}
($i$) The existence of the long exact sequence is a direct consequence of the Zig-Zag Lemma of homological algebra (see, e.g., \cite{Weibel94}), applied to the short exact sequence of Koszul complexes \eqref{eq:koszul_chains}; the exactness of the rows and the commutativity of the squares are guaranteed by the construction of $i_k, p_k$ and by Lemma \ref{lem:comm_diag}. For the reader's convenience and to clarify the underlying structure, we provide an explicit description of the connecting homomorphism $\partial_k \colon H^k(\mathbf{B}) \to H^{k+1}(\mathbf{A})$.
Let $[\eta] \in H^k(\mathbf{B})$ be represented by a cocycle $\eta \in \ker \delta_{\mathbf{B}}^k$. To define $\partial_k([\eta])$, we first select the most natural pre-image of $\eta$ under the projection $p_k$, namely $\zeta = \begin{pmatrix} 0 \\ \eta \end{pmatrix} \in \Lambda^k[\mathcal{X} \oplus \mathcal{Y}]$. Applying the boundary map $\delta_{\mathbf{M}}^k$ to this vector, we obtain
\[
\delta_{\mathbf{M}}^k \begin{pmatrix} 0 \\ \eta \end{pmatrix} = \begin{pmatrix} \delta_{\mathbf{A}}^k & \delta_{\mathbf{C}}^k \\ 0 & \delta_{\mathbf{B}}^k \end{pmatrix} \begin{pmatrix} 0 \\ \eta \end{pmatrix} = \begin{pmatrix} \delta_{\mathbf{C}}^k \eta \\ \delta_{\mathbf{B}}^k \eta \end{pmatrix} = \begin{pmatrix} \delta_{\mathbf{C}}^k \eta \\ 0 \end{pmatrix},
\]
where the second component vanishes because $\delta_{\mathbf{B}}^k \eta = 0$. Since $\begin{pmatrix} \delta_{\mathbf{C}}^k \eta \\ 0 \end{pmatrix} = i_{k+1}(\delta_{\mathbf{C}}^k \eta)$, this vector corresponds to $\delta_{\mathbf{C}}^k \eta \in \Lambda^{k+1}[\mathcal{X}]$. We then define the connecting homomorphism by mapping the class $[\eta]$ to the cohomology class of $\delta_{\mathbf{C}}^k \eta$:
\begin{equation}\label{eq:partial_def_multivar}
    \partial_k([\eta]) = [\delta_{\mathbf{C}}^k \eta] \in H^{k+1}(\mathbf{A}).
\end{equation}
This shows that the connecting homomorphism is precisely the map $\partial_k$ of Theorem \ref{thm:M_range_closed}, induced by the coupling tuple $\mathbf{C}$ on the level of cohomology.
To verify that this definition is independent of the choice of pre-image, let $\zeta' = \begin{pmatrix} \xi \\ \eta \end{pmatrix}$ be any other pre-image of $\eta$ under $p_k$, for some $\xi \in \Lambda^k[\mathcal{X}]$. Applying $\delta_{\mathbf{M}}^k$ yields
\[
\delta_{\mathbf{M}}^k \begin{pmatrix} \xi \\ \eta \end{pmatrix} = \begin{pmatrix} \delta_{\mathbf{A}}^k \xi + \delta_{\mathbf{C}}^k \eta \\ \delta_{\mathbf{B}}^k \eta \end{pmatrix} = \begin{pmatrix} \delta_{\mathbf{A}}^k \xi + \delta_{\mathbf{C}}^k \eta \\ 0 \end{pmatrix} = i_{k+1}(\delta_{\mathbf{A}}^k \xi + \delta_{\mathbf{C}}^k \eta).
\]
Taking the cohomology class in $H^{k+1}(\mathbf{A})$, the term $\delta_{\mathbf{A}}^k \xi$ vanishes. Thus,
\[
[\delta_{\mathbf{A}}^k \xi + \delta_{\mathbf{C}}^k \eta] = [\delta_{\mathbf{A}}^k \xi] + [\delta_{\mathbf{C}}^k \eta] = 0 + [\delta_{\mathbf{C}}^k \eta] = [\delta_{\mathbf{C}}^k \eta],
\]
which confirms that the definition of $\partial_k([\eta])$ is independent of $\xi$.
The fact that $\partial_k$ is independent of the choice of representative $\eta$ within the class $[\eta]$ is already shown in the proof of Theorem \ref{thm:M_range_closed}, and thus the map $\partial_k$ is well-defined.

($ii$) Under the closedness assumptions, all the quotients $H^k(\mathbf{A})$, $H^k(\mathbf{M})$, $H^k(\mathbf{B})$ are Banach spaces with respect to the quotient norms. The maps $\iota_k^*$ and $\pi_k^*$ are induced by the bounded operators $i_k$ and $p_k$ on the corresponding subquotients, and are therefore bounded. As for the connecting homomorphism, if $\eta' = \eta + \delta_{\mathbf{B}}^{k-1}\zeta$ is an arbitrary representative of the class $[\eta] \in H^k(\mathbf{B})$, then $[\delta_{\mathbf{C}}^k \eta] = [\delta_{\mathbf{C}}^k \eta']$ by the well-definedness argument above, whence
\[
\|\partial_k[\eta]\| = \|[\delta_{\mathbf{C}}^k \eta']\| \leq \|\delta_{\mathbf{C}}^k \eta'\| \leq \|\delta_{\mathbf{C}}^k\|\,\|\eta'\|.
\]
Taking the infimum over all representatives $\eta'$ gives $\|\partial_k[\eta]\| \leq \|\delta_{\mathbf{C}}^k\|\,\|[\eta]\|$.
\end{proof}
We now turn to the characterization of the Taylor invertibility of $\mathbf{M}$ in terms of the connecting homomorphisms. Closed-range hypotheses enter only in its second part, which upgrades the linear bijections to topological isomorphisms.
\begin{theorem}\label{thm:taylor_char}
Let $\mathcal{X}$ and $\mathcal{Y}$ be Banach spaces, and let $\mathbf{M}$ be a commuting $n$-tuple of upper triangular operator matrices acting on $\mathcal{X} \oplus \mathcal{Y}$ as defined in \eqref{eq:M_tuple}. Then the $n$-tuple $\mathbf{M}$ is Taylor invertible if and only if the map
\[
\partial_k: H^k(\mathbf{B}) \to H^{k+1}(\mathbf{A}), \qquad \partial_k([\eta]) = [\delta_{\mathbf{C}}^k \eta], \quad \eta \in \ker \delta_{\mathbf{B}}^k,
\]
is a linear bijection for every $k\in\mathbb{Z}$.

If, in addition, the ranges $\mathcal{R}(\delta_{\mathbf{A}}^k)$ and $\mathcal{R}(\delta_{\mathbf{B}}^k)$ are closed in $\Lambda^{k+1}[\mathcal{X}]$ and $\Lambda^{k+1}[\mathcal{Y}]$, respectively, for every $k \in \mathbb{Z}$, then $\mathbf{M}$ is Taylor invertible if and only if each $\partial_k$ is a topological isomorphism of Banach spaces.
\end{theorem}
\begin{proof}
$(\implies)$ Suppose that the commuting $n$-tuple $\mathbf{M}$ is Taylor invertible, so that $H^k(\mathbf{M}) = 0$ for all $k \in \mathbb{Z}$. The long exact sequence \eqref{eq:zig_zag_Taylor} of Lemma \ref{lem:zig_zag_Taylor}($i$) then contains, for each $k$, the exact segment
\begin{equation*}
    \begin{tikzcd}[column sep=2.2em]
        0 = H^k(\mathbf{M}) \arrow[r, "\pi_k^*"] & H^k(\mathbf{B}) \arrow[r, "\partial_k"]  &
        H^{k+1}(\mathbf{A}) \arrow[r, "\iota_{k+1}^*"]  & H^{k+1}(\mathbf{M}) = 0.
    \end{tikzcd}
\end{equation*}
Exactness at $H^k(\mathbf{B})$ gives $\mathcal{N}(\partial_k) = \mathcal{R}(\pi_k^*) = \{0\}$, so $\partial_k$ is injective; exactness at $H^{k+1}(\mathbf{A})$ gives $\mathcal{R}(\partial_k) = \mathcal{N}(\iota_{k+1}^*) = H^{k+1}(\mathbf{A})$, so $\partial_k$ is surjective. Hence each $\partial_k$ is a linear bijection.
\medskip

$(\impliedby)$ Conversely, assume that $\partial_k$ is a linear bijection for every $k \in \mathbb{Z}$. To establish the Taylor invertibility of $\mathbf{M}$, it suffices to show that its cohomology groups vanish for all $k$. Consider the following segment of \eqref{eq:zig_zag_Taylor}:
\begin{equation*}
    \begin{tikzcd}[column sep=2.5em]
        H^{k-1}(\mathbf{B}) \arrow[r, "\partial_{k-1}"] & H^k(\mathbf{A}) \arrow[r, "\iota_k^*"] & H^k(\mathbf{M}) \arrow[r, "\pi_k^*"] & H^k(\mathbf{B}) \arrow[r, "\partial_k"] & H^{k+1}(\mathbf{A})
    \end{tikzcd}
\end{equation*}
The connecting homomorphism $\partial_{k-1} \colon H^{k-1}(\mathbf{B}) \to H^k(\mathbf{A})$ is a bijection. In particular, it is surjective, so by the exactness of the sequence at $H^k(\mathbf{A})$, we have
\[
\mathcal{N}(\iota_k^*) = \mathcal{R}(\partial_{k-1}) = H^k(\mathbf{A}).
\]
This implies that $\iota_k^* = 0$, and thus $\mathcal{R}(\iota_k^*) = \{0\}$.
On the other hand, since $\partial_k \colon H^k(\mathbf{B}) \to H^{k+1}(\mathbf{A})$ is also a bijection, it is injective. By exactness at $H^k(\mathbf{B})$, we obtain
\[
 \mathcal{R}(\pi_k^*)=\mathcal{N}(\partial_k) = \{0\},
\]
which implies that $\pi_k^* = 0$, and thus $\mathcal{N}(\pi_k^*) = H^k(\mathbf{M})$. Combining these observations with the exactness of the sequence at $H^k(\mathbf{M})$, it follows that
\begin{equation*}
    H^k(\mathbf{M}) = \mathcal{N}(\pi_k^*) = \mathcal{R}(\iota_k^*) = \{0\}.
\end{equation*}
Thus, we conclude that $\mathbf{M}$ is Taylor invertible.
\medskip

Finally, assume that the ranges $\mathcal{R}(\delta_{\mathbf{A}}^k)$ and $\mathcal{R}(\delta_{\mathbf{B}}^k)$ are closed for every $k \in \mathbb{Z}$. If $\mathbf{M}$ is Taylor invertible, then the ranges $\mathcal{R}(\delta_{\mathbf{M}}^k) = \ker \delta_{\mathbf{M}}^{k+1}$ are closed as well, so by Lemma \ref{lem:zig_zag_Taylor}($ii$) all the cohomology spaces are Banach and each $\partial_k$ is bounded; being bijective by the first part, $\partial_k$ is a topological isomorphism by the Open Mapping Theorem. Conversely, if each $\partial_k$ is a topological isomorphism, it is in particular a linear bijection, and the first part of the theorem yields the Taylor invertibility of $\mathbf{M}$.
\end{proof}
\begin{remark}\label{rem:edge_degrees}
    With the convention of Section \ref{sec:prelim}, the bijectivity of the two extreme connecting maps in Theorem \ref{thm:taylor_char} carries the following content: $\partial_{-1} \colon H^{-1}(\mathbf{B}) = \{0\} \to H^0(\mathbf{A})$ is bijective if and only if $H^0(\mathbf{A}) = \{0\}$, i.e., $\delta_{\mathbf{A}}^0$ is injective, while $\partial_n \colon H^n(\mathbf{B}) \to H^{n+1}(\mathbf{A}) = \{0\}$ is bijective if and only if $H^n(\mathbf{B}) = \{0\}$, i.e., $\delta_{\mathbf{B}}^{n-1}$ is surjective. Thus Theorem \ref{thm:taylor_char} subsumes parts ($i$) and ($iii$) of Theorem \ref{thm:AB_range_closed}; the additional content of the latter is its topological part ($ii$), namely that $\delta_{\mathbf{A}}^0$ is bounded below and not merely injective.
\end{remark}
\begin{remark}\label{rem:bijective_enough}
    Let us record the relation between Theorem \ref{thm:taylor_char} and Theorem \ref{thm:M_range_closed}. Once $\mathbf{M}$ is Taylor invertible, the ranges $\mathcal{R}(\delta_{\mathbf{M}}^k) = \ker \delta_{\mathbf{M}}^{k+1}$ are closed for trivial reasons; the interest of Theorem \ref{thm:M_range_closed} is that it operates at a single degree $k$ and, under the closedness of the ranges of $\delta_{\mathbf{A}}^k$ and $\delta_{\mathbf{B}}^k$, characterizes the closedness of $\mathcal{R}(\delta_{\mathbf{M}}^k)$ by that of $\mathcal{R}(\partial_k)$, without any exactness or bijectivity assumptions.
\end{remark}
The long exact sequence also yields the classical two-out-of-three principle for the triple $(\mathbf{A}, \mathbf{B}, \mathbf{M})$, together with the accompanying spectral identity. The principle goes back to Taylor \cite{Taylor70}, who derived from the long exact sequence that the union of any two of the sets $\sigma_T(\mathbf{M})$, $\sigma_T(\mathbf{A})$, $\sigma_T(\mathbf{B})$ contains the third, and the two implications in which $\mathbf{M}$, respectively $\mathbf{B}$, is the conclusion were also obtained by Curto \cite{Curto80} in the Hilbert space setting. We record the short derivation from \eqref{eq:zig_zag_Taylor}, which uses only Lemma \ref{lem:zig_zag_Taylor}($i$).
\begin{corollary}\label{cor:two_out_of_three}
    Let $\mathbf{A} \in \mathfrak{B}(\mathcal{X})^n$, $\mathbf{B} \in \mathfrak{B}(\mathcal{Y})^n$, and $\mathbf{C} \in \mathfrak{B}(\mathcal{Y}, \mathcal{X})^n$ be such that $\mathbf{M} = \begin{pmatrix} \mathbf{A} & \mathbf{C} \\ \mathbf{0} & \mathbf{B} \end{pmatrix}$ is a commuting $n$-tuple.
    \begin{enumerate}[\rm($i$)]
        \item If any two of the tuples $\mathbf{A}$, $\mathbf{B}$, $\mathbf{M}$ are Taylor invertible, then so is the third.
        \item Each of the sets $\sigma_T(\mathbf{A})$, $\sigma_T(\mathbf{B})$, $\sigma_T(\mathbf{M})$ is contained in the union of the other two. Consequently,
        \begin{equation}\label{eq:spectral_identity}
            \sigma_T(\mathbf{A}) \cup \sigma_T(\mathbf{B}) = \sigma_T(\mathbf{M}) \cup \bigl(\sigma_T(\mathbf{A}) \cap \sigma_T(\mathbf{B})\bigr).
        \end{equation}
    \end{enumerate}
\end{corollary}
\begin{proof}
    ($i$) Fix $k \in \mathbb{Z}$ and consider the exact segment
    \begin{equation*}
        \begin{tikzcd}[column sep=2.5em]
            H^{k-1}(\mathbf{B}) \arrow[r, "\partial_{k-1}"] & H^k(\mathbf{A}) \arrow[r, "\iota_k^*"] & H^k(\mathbf{M}) \arrow[r, "\pi_k^*"] & H^k(\mathbf{B}) \arrow[r, "\partial_k"] & H^{k+1}(\mathbf{A})
        \end{tikzcd}
    \end{equation*}
    of the long exact sequence from Lemma \ref{lem:zig_zag_Taylor}($i$). If $\mathbf{A}$ and $\mathbf{B}$ are Taylor invertible, then $H^k(\mathbf{A}) = \{0\}$ and $H^k(\mathbf{B}) = \{0\}$, so exactness at $H^k(\mathbf{M})$ gives $H^k(\mathbf{M}) = \mathcal{N}(\pi_k^*) = \mathcal{R}(\iota_k^*) = \{0\}$. If $\mathbf{A}$ and $\mathbf{M}$ are Taylor invertible, then exactness at $H^k(\mathbf{B})$ gives $\mathcal{N}(\partial_k) = \mathcal{R}(\pi_k^*) = \{0\}$, and since $\partial_k$ maps into $H^{k+1}(\mathbf{A}) = \{0\}$, we obtain $H^k(\mathbf{B}) = \mathcal{N}(\partial_k) = \{0\}$. If $\mathbf{B}$ and $\mathbf{M}$ are Taylor invertible, then exactness at $H^k(\mathbf{A})$ gives $\mathcal{N}(\iota_k^*) = \mathcal{R}(\partial_{k-1}) = \{0\}$, and since $\iota_k^*$ maps into $H^k(\mathbf{M}) = \{0\}$, we obtain $H^k(\mathbf{A}) = \mathcal{N}(\iota_k^*) = \{0\}$. Since $k$ was arbitrary, in each case the third tuple is Taylor invertible.

    ($ii$) For every $\lambda = (\lambda_1, \dots, \lambda_n) \in \mathbb{C}^n$, the translated tuple $\mathbf{M} - \lambda$ is again a commuting $n$-tuple of the form \eqref{eq:M_tuple}, with diagonal tuples $\mathbf{A} - \lambda$ and $\mathbf{B} - \lambda$ and the same coupling tuple $\mathbf{C}$; indeed, the terms $-\lambda_i C_j + \lambda_j C_i$ appear on both sides of \eqref{eq:coupling_relation} after the translation. The three inclusions are therefore the contrapositives of the three cases of ($i$), applied to $\mathbf{M} - \lambda$. In \eqref{eq:spectral_identity}, the inclusion $\supseteq$ follows from $\sigma_T(\mathbf{M}) \subseteq \sigma_T(\mathbf{A}) \cup \sigma_T(\mathbf{B})$ together with the trivial inclusion $\sigma_T(\mathbf{A}) \cap \sigma_T(\mathbf{B}) \subseteq \sigma_T(\mathbf{A}) \cup \sigma_T(\mathbf{B})$. For $\subseteq$, let $\lambda \in \sigma_T(\mathbf{A}) \setminus \sigma_T(\mathbf{M})$. Since $\sigma_T(\mathbf{A}) \subseteq \sigma_T(\mathbf{M}) \cup \sigma_T(\mathbf{B})$, we get $\lambda \in \sigma_T(\mathbf{B})$, hence $\lambda \in \sigma_T(\mathbf{A}) \cap \sigma_T(\mathbf{B})$. The case $\lambda \in \sigma_T(\mathbf{B}) \setminus \sigma_T(\mathbf{M})$ is symmetric.
\end{proof}
\bigskip
\section{The Multivariable Completion Problem}\label{sec:completion}
\bigskip
Having characterized in Theorem \ref{thm:taylor_char} the Taylor invertibility of a \emph{given} commuting matrix tuple, we now turn to the completion problem itself: starting from the diagonal tuples $\mathbf{A}$ and $\mathbf{B}$, we ask whether the partially specified tuple $\begin{pmatrix} \mathbf{A} & \mathbf{?} \\ \mathbf{0} & \mathbf{B} \end{pmatrix}$ admits a coupling $\mathbf{C}$ rendering $\mathbf{M}_{\mathbf{C}}$ a commuting, Taylor invertible $n$-tuple. Since the coupling tuple now varies, we write $\mathbf{M}_{\mathbf{C}}$ for the tuple \eqref{eq:M_tuple} and $\partial_k^{\mathbf{C}}$ for the maps of Theorem \ref{thm:taylor_char}, in order to emphasize their dependence on $\mathbf{C}$. This is the multivariable analogue of Theorem \ref{thm:1D_completion}. As in Theorem \ref{thm:taylor_char}, closed-range hypotheses enter only when the criterion is to be phrased in terms of topological isomorphisms.
\begin{theorem}\label{thm:completion}
    Let $\mathbf{A} \in \mathfrak{B}(\mathcal{X})^n$ and $\mathbf{B} \in \mathfrak{B}(\mathcal{Y})^n$ be commuting $n$-tuples. Then the following statements are equivalent:
    \begin{enumerate}[\rm($i$)]
        \item There exists an $n$-tuple $\mathbf{C} = (C_1, \dots, C_n)$ of bounded linear operators such that
        \[
        \mathbf{M}_{\mathbf{C}} = \begin{pmatrix} \mathbf{A} & \mathbf{C} \\ \mathbf{0} & \mathbf{B} \end{pmatrix}
        \]
        is a commuting and Taylor invertible $n$-tuple.
        \item There exists an $n$-tuple $\mathbf{C} = (C_1, \dots, C_n)$ of bounded linear operators satisfying the intertwining relations
        \begin{equation}\label{eq:completion_coupling}
            A_i C_j - A_j C_i = C_j B_i - C_i B_j \qquad \text{for all } i, j,
        \end{equation}
        whose induced maps on cohomology
        \[
        \partial_k^{\mathbf{C}} \colon H^k(\mathbf{B}) \to H^{k+1}(\mathbf{A}), \qquad \partial_k^{\mathbf{C}}([\eta]) = [\delta_{\mathbf{C}}^k \eta],
        \]
        are linear bijections for every $k \in \mathbb{Z}$.
    \end{enumerate}
    Moreover, when these conditions hold, one and the same coupling tuple $\mathbf{C}$ realizes both ($i$) and ($ii$). If, in addition, the ranges $\mathcal{R}(\delta_{\mathbf{A}}^k)$ and $\mathcal{R}(\delta_{\mathbf{B}}^k)$ are closed for every $k \in \mathbb{Z}$, then the maps $\partial_k^{\mathbf{C}}$ in ($ii$) may equivalently be required to be topological isomorphisms of Banach spaces.
\end{theorem}
\begin{proof}
    As recorded in \eqref{eq:coupling_relation}, the tuple $\mathbf{M}_{\mathbf{C}}$ is commuting if and only if $\mathbf{C}$ satisfies \eqref{eq:completion_coupling}. For every such $\mathbf{C}$, the first part of Theorem \ref{thm:taylor_char} shows that $\mathbf{M}_{\mathbf{C}}$ is Taylor invertible if and only if all the maps $\partial_k^{\mathbf{C}}$ are linear bijections, and, under the closed-range hypotheses, the second part shows that this is the case if and only if all of them are topological isomorphisms of Banach spaces. This proves the equivalence of ($i$) and ($ii$), with the same coupling tuple $\mathbf{C}$ in both, as well as the last assertion.
\end{proof}
The connecting homomorphisms carry more structure than that of mere linear maps, and this additional structure yields a finer necessary condition for the existence of a Taylor invertible completion. For $1 \leq i \leq n$ and a commuting $n$-tuple $\mathbf{T}$, the creation operator $E_i$ maps $\Lambda^k[\mathcal{X}]$ into $\Lambda^{k+1}[\mathcal{X}]$ and anticommutes with $\delta_{\mathbf{T}}$, since
\[
\delta_{\mathbf{T}} E_i (x_F \otimes e_F) = \sum_{j=1}^n T_j x_F \otimes (e_j \wedge e_i \wedge e_F) = - \sum_{j=1}^n T_j x_F \otimes (e_i \wedge e_j \wedge e_F) = - E_i \delta_{\mathbf{T}} (x_F \otimes e_F).
\]
Hence $E_i$ maps $\ker \delta_{\mathbf{T}}^k$ into $\ker \delta_{\mathbf{T}}^{k+1}$ and $\mathcal{R}(\delta_{\mathbf{T}}^{k-1})$ into $\mathcal{R}(\delta_{\mathbf{T}}^{k})$, and therefore induces a linear map
\[
\varepsilon_i^{\mathbf{T}} \colon H^k(\mathbf{T}) \to H^{k+1}(\mathbf{T}), \qquad \varepsilon_i^{\mathbf{T}} [\xi] = [E_i \xi].
\]
Since $E_i E_j = - E_j E_i$ for all $i, j$, the maps $\varepsilon_1^{\mathbf{T}}, \dots, \varepsilon_n^{\mathbf{T}}$ make $H^\bullet(\mathbf{T}) = \bigoplus_{k \in \mathbb{Z}} H^k(\mathbf{T})$ a graded module over the exterior algebra $\Lambda_n[\mathbf{e}]$, the generator $e_i$ acting as $\varepsilon_i^{\mathbf{T}}$. We keep the superscript, since the maps $\varepsilon_i^{\mathbf{A}}$ and $\varepsilon_i^{\mathbf{B}}$, which act on the cohomology of different tuples, will appear in the same identity. For a commuting $n$-tuple $\mathbf{A}$ we write $H^{\bullet+1}(\mathbf{A})$ for the graded $\Lambda_n[\mathbf{e}]$-module whose component of degree $k$ is $H^{k+1}(\mathbf{A})$, the generator $e_i$ acting as $\varepsilon_i^{\mathbf{A}}$. The same computation applies verbatim to $\delta_{\mathbf{C}} = \sum_{j=1}^n C_j \otimes E_j$ and gives
\[
\delta_{\mathbf{C}}^{k+1} E_i = - E_i \delta_{\mathbf{C}}^k \quad \text{on } \Lambda^k[\mathcal{Y}],
\]
where $E_i$ acts on $\Lambda^k[\mathcal{Y}]$ on the left-hand side and on $\Lambda^{k+1}[\mathcal{X}]$ on the right-hand side.
\begin{theorem}\label{thm:module}
    Let $\mathbf{M}$ be the commuting $n$-tuple defined by \eqref{eq:M_tuple}. Then, for all $k \in \mathbb{Z}$ and $1 \leq i \leq n$,
    \[
    \partial_{k+1} \circ \varepsilon_i^{\mathbf{B}} = - \varepsilon_i^{\mathbf{A}} \circ \partial_k
    \]
    as linear maps from $H^k(\mathbf{B})$ to $H^{k+2}(\mathbf{A})$.
    Consequently, if $\mathbf{M}$ is Taylor invertible, then the maps $(-1)^k \partial_k$, $k \in \mathbb{Z}$, constitute an isomorphism of graded $\Lambda_n[\mathbf{e}]$-modules from $H^\bullet(\mathbf{B})$ onto $H^{\bullet+1}(\mathbf{A})$. In particular, for all $k \in \mathbb{Z}$ and $1 \leq i \leq n$, the map $\varepsilon_i^{\mathbf{B}} \colon H^k(\mathbf{B}) \to H^{k+1}(\mathbf{B})$ is injective, respectively surjective, respectively zero, if and only if the map $\varepsilon_i^{\mathbf{A}} \colon H^{k+1}(\mathbf{A}) \to H^{k+2}(\mathbf{A})$ has the same property.
\end{theorem}
\begin{proof}
    Let $\eta \in \ker \delta_{\mathbf{B}}^k$. Then $E_i \eta \in \ker \delta_{\mathbf{B}}^{k+1}$, and the anticommutation relation for $\delta_{\mathbf{C}}$ gives
    \[
    \partial_{k+1}(\varepsilon_i^{\mathbf{B}} [\eta]) = \partial_{k+1} [E_i \eta] = [\delta_{\mathbf{C}}^{k+1} E_i \eta] = [- E_i \delta_{\mathbf{C}}^k \eta] = - \varepsilon_i^{\mathbf{A}} [\delta_{\mathbf{C}}^k \eta] = - \varepsilon_i^{\mathbf{A}} \partial_k [\eta].
    \]
    Now assume that $\mathbf{M}$ is Taylor invertible, so that each $\partial_k$ is a linear bijection by Theorem \ref{thm:taylor_char}, and put $\partial_k' = (-1)^k \partial_k$. Then
    \[
    \partial_{k+1}' \circ \varepsilon_i^{\mathbf{B}} = (-1)^{k+1} \partial_{k+1} \circ \varepsilon_i^{\mathbf{B}} = (-1)^{k+1} (- \varepsilon_i^{\mathbf{A}} \circ \partial_k) = \varepsilon_i^{\mathbf{A}} \circ \partial_k',
    \]
    so the family $(\partial_k')_{k \in \mathbb{Z}}$ commutes with the action of every generator $e_i$ and is therefore a homomorphism of graded $\Lambda_n[\mathbf{e}]$-modules from $H^\bullet(\mathbf{B})$ to $H^{\bullet+1}(\mathbf{A})$. Being bijective in every degree, it is an isomorphism of graded modules, its inverse being automatically $\Lambda_n[\mathbf{e}]$-linear. Finally, the identity $\varepsilon_i^{\mathbf{B}} = (\partial_{k+1}')^{-1} \circ \varepsilon_i^{\mathbf{A}} \circ \partial_k'$ on $H^k(\mathbf{B})$ shows that the two maps in the last assertion are simultaneously injective, surjective or zero.
\end{proof}
Theorem \ref{thm:module} is purely algebraic. Under the closed-range hypotheses of the second part of Theorem \ref{thm:taylor_char}, the maps $(-1)^k \partial_k$ are moreover topological isomorphisms, and the maps $\varepsilon_i^{\mathbf{A}}$ and $\varepsilon_i^{\mathbf{B}}$ are bounded, being induced by the bounded operators $E_i$ on the corresponding subquotients.

Theorems \ref{thm:AB_range_closed} and \ref{thm:module} yield necessary conditions for the existence of a Taylor invertible completion which involve only the diagonal tuples.
\begin{corollary}\label{cor:necessary}
    Let $\mathbf{A} \in \mathfrak{B}(\mathcal{X})^n$ and $\mathbf{B} \in \mathfrak{B}(\mathcal{Y})^n$ be commuting $n$-tuples, and suppose that there exists a coupling tuple $\mathbf{C}$ for which $\mathbf{M}_{\mathbf{C}}$ is a commuting and Taylor invertible $n$-tuple. Then:
    \begin{enumerate}[\rm($i$)]
        \item $\delta_{\mathbf{A}}^0$ is bounded below;
        \item $H^\bullet(\mathbf{B})$ and $H^{\bullet+1}(\mathbf{A})$ are isomorphic as graded $\Lambda_n[\mathbf{e}]$-modules; in particular, $H^k(\mathbf{B})$ and $H^{k+1}(\mathbf{A})$ are linearly isomorphic for every $k \in \mathbb{Z}$.
    \end{enumerate}
    In particular, $\delta_{\mathbf{B}}^{n-1}$ is surjective.
\end{corollary}
\begin{proof}
    Condition ($i$) follows from Theorem \ref{thm:AB_range_closed} applied to $\mathbf{M}_{\mathbf{C}}$, and ($ii$) follows from Theorem \ref{thm:module}. For $k = n$, condition ($ii$) gives $H^n(\mathbf{B}) \cong H^{n+1}(\mathbf{A}) = \{0\}$, that is, $\mathcal{R}(\delta_{\mathbf{B}}^{n-1}) = \Lambda^n[\mathcal{Y}]$.
\end{proof}
As in Remark \ref{rem:edge_degrees}, the case $k = -1$ of ($ii$) gives $H^0(\mathbf{A}) = \{0\}$, that is, the injectivity of $\delta_{\mathbf{A}}^0$. Thus the content of ($i$) beyond ($ii$) is the closedness of $\mathcal{R}(\delta_{\mathbf{A}}^0)$. At the top degree no such addition is needed, since the range of the surjective map $\delta_{\mathbf{B}}^{n-1}$ is closed.

For $n = 1$, condition ($ii$) says that $\mathcal{N}(A) = H^0(\mathbf{A})$ and $\mathcal{Y}/\mathcal{R}(B) = H^1(\mathbf{B})$ vanish and that $\mathcal{N}(B) = H^0(\mathbf{B})$ and $\mathcal{X}/\mathcal{R}(A) = H^1(\mathbf{A})$ are linearly isomorphic, the generator $e_1$ necessarily acting as zero on these spaces. Together with ($i$), this means that $A$ is bounded below, $B$ is surjective and $\mathcal{N}(B) \cong \mathcal{X}/\mathcal{R}(A)$ as vector spaces. Thus Corollary \ref{cor:necessary} is the multivariable counterpart of the three conditions of Theorem \ref{thm:1D_completion}, and the module structure in ($ii$) is a genuinely multivariable feature. Let us stress that the isomorphisms in ($ii$) are linear isomorphisms. Theorem \ref{thm:1D_completion} shows that for $n = 1$ these conditions become sufficient when the boundedness from below of $A$ and the surjectivity of $B$ are strengthened to the left invertibility of $A$ and the right invertibility of $B$, respectively, and the isomorphism $\mathcal{N}(B) \cong \mathcal{X}/\mathcal{R}(A)$ is required to be an isomorphism of Banach spaces. The distinction is essential: the closedness of the ranges alone does not turn an abstract linear isomorphism into a topological one, although, by the second part of Theorem \ref{thm:taylor_char}, the particular isomorphisms $\partial_k^{\mathbf{C}}$ arising from a Taylor invertible completion are topological whenever the relevant ranges are closed. Example \ref{ex:no_completion} below shows that no such strengthening restores sufficiency for $n = 2$. Before comparing the two situations further, we record the freedom one has in modifying the coupling tuple.
\begin{proposition}\label{prop:gauge}
    Let $\mathbf{A} \in \mathfrak{B}(\mathcal{X})^n$ and $\mathbf{B} \in \mathfrak{B}(\mathcal{Y})^n$ be commuting $n$-tuples, let $\mathbf{C} \in \mathfrak{B}(\mathcal{Y}, \mathcal{X})^n$ satisfy \eqref{eq:completion_coupling}, and let $D \in \mathfrak{B}(\mathcal{Y}, \mathcal{X})$. Put $\mathbf{C}' = (C_1', \dots, C_n')$, where $C_i' = C_i + A_i D - D B_i$ for $1 \leq i \leq n$. Then:
    \begin{enumerate}[\rm($i$)]
        \item $\mathbf{C}'$ satisfies \eqref{eq:completion_coupling};
        \item $\mathbf{M}_{\mathbf{C}'} = S^{-1} \mathbf{M}_{\mathbf{C}} S$, where $S = \begin{pmatrix} I & D \\ 0 & I \end{pmatrix} \in \mathfrak{B}(\mathcal{X} \oplus \mathcal{Y})$;
        \item $\partial_k^{\mathbf{C}'} = \partial_k^{\mathbf{C}}$ for every $k \in \mathbb{Z}$.
    \end{enumerate}
    In particular, $\mathbf{M}_{\mathbf{C}'}$ is Taylor invertible if and only if $\mathbf{M}_{\mathbf{C}}$ is Taylor invertible.
\end{proposition}
\begin{proof}
    ($ii$) The operator $S$ is invertible with $S^{-1} = \begin{pmatrix} I & -D \\ 0 & I \end{pmatrix}$, and for every $i$ we have
    \begin{align*}
        S^{-1} M_i S &= \begin{pmatrix} I & -D \\ 0 & I \end{pmatrix} \begin{pmatrix} A_i & C_i \\ 0 & B_i \end{pmatrix} \begin{pmatrix} I & D \\ 0 & I \end{pmatrix} = \begin{pmatrix} A_i & C_i - D B_i \\ 0 & B_i \end{pmatrix} \begin{pmatrix} I & D \\ 0 & I \end{pmatrix} \\
        &= \begin{pmatrix} A_i & A_i D + C_i - D B_i \\ 0 & B_i \end{pmatrix} = \begin{pmatrix} A_i & C_i' \\ 0 & B_i \end{pmatrix}.
    \end{align*}

    ($i$) By ($ii$), the tuple $\mathbf{M}_{\mathbf{C}'}$ is similar to the commuting $n$-tuple $\mathbf{M}_{\mathbf{C}}$, and is therefore itself a commuting $n$-tuple. As recorded in \eqref{eq:coupling_relation}, this means that $\mathbf{C}'$ satisfies \eqref{eq:completion_coupling}.

    ($iii$) Extending $D$ diagonally to the exterior algebra as in Section \ref{sec:prelim}, for $\eta = \sum_{|F|=k} y_F \otimes e_F \in \Lambda^k[\mathcal{Y}]$ we have
    \[
    \delta_{\mathbf{C}'}^k \eta = \sum_{F} \sum_{i=1}^n (C_i + A_i D - D B_i) y_F \otimes (e_i \wedge e_F) = \delta_{\mathbf{C}}^k \eta + \delta_{\mathbf{A}}^k (D \eta) - D (\delta_{\mathbf{B}}^k \eta).
    \]
    If $\eta \in \ker \delta_{\mathbf{B}}^k$, this reduces to $\delta_{\mathbf{C}'}^k \eta = \delta_{\mathbf{C}}^k \eta + \delta_{\mathbf{A}}^k (D \eta)$, and since $\delta_{\mathbf{A}}^k (D \eta) \in \mathcal{R}(\delta_{\mathbf{A}}^k)$, we obtain $[\delta_{\mathbf{C}'}^k \eta] = [\delta_{\mathbf{C}}^k \eta]$ in $H^{k+1}(\mathbf{A})$, that is, $\partial_k^{\mathbf{C}'}([\eta]) = \partial_k^{\mathbf{C}}([\eta])$.

    The final assertion follows from ($iii$) and the first part of Theorem \ref{thm:taylor_char}.
\end{proof}
In other words, both the Taylor invertibility of $\mathbf{M}_{\mathbf{C}}$ and the connecting homomorphisms depend only on the class of $\mathbf{C}$ modulo the couplings of the form $(A_i D - D B_i)_{i=1}^n$ with $D \in \mathfrak{B}(\mathcal{Y}, \mathcal{X})$, which we call trivial. Since $\partial_k^{\mathbf{0}} = 0$ for all $k$, Theorem \ref{thm:taylor_char} shows that $\mathbf{M}_{\mathbf{0}} = \mathbf{A} \oplus \mathbf{B}$, and hence $\mathbf{M}_{\mathbf{C}}$ for every trivial coupling $\mathbf{C}$, is Taylor invertible if and only if $\mathbf{A}$ and $\mathbf{B}$ are both Taylor invertible.
The following example shows that condition ($ii$) of Corollary \ref{cor:necessary} cannot be weakened to the existence of isomorphisms $H^k(\mathbf{B}) \cong H^{k+1}(\mathbf{A})$ of Banach spaces, and that the latter, together with ($i$), do not guarantee the existence of a Taylor invertible completion.
\begin{example}\label{ex:no_completion}
    Let $n = 2$, let $\mathcal{X} = \mathcal{Y} = H^2(\mathbb{D})$ be the Hardy space of the unit disc, and let $S$ be the unilateral shift, $(Sf)(z) = z f(z)$. Recall that $S$ is an isometry, that $S^*$ is surjective with $\ker S^* = \mathbb{C} 1$, where $1$ denotes the constant function, and that $\mathcal{R}(S) = \{ f \in H^2(\mathbb{D}) : f(0) = 0 \}$ has codimension one. Consider the commuting pairs
    \[
    \mathbf{A} = (S, 0) \qquad \text{and} \qquad \mathbf{B} = (0, S^*).
    \]
    For $\mathbf{A}$ we have $\delta_{\mathbf{A}}^0 x = Sx \otimes e_1$ and $\delta_{\mathbf{A}}^1 (x_1 \otimes e_1 + x_2 \otimes e_2) = S x_2 \otimes (e_1 \wedge e_2)$. Hence $\delta_{\mathbf{A}}^0$ is bounded below, $\ker \delta_{\mathbf{A}}^1 = \{ x \otimes e_1 : x \in \mathcal{X} \}$, and
    \begin{align*}
    H^0(\mathbf{A}) &= \{0\}, \\
    H^1(\mathbf{A}) &= \{ x \otimes e_1 : x \in \mathcal{X} \} / \{ Sx \otimes e_1 : x \in \mathcal{X} \} \cong \mathcal{X}/\mathcal{R}(S), \\
    H^2(\mathbf{A}) &= \Lambda^2[\mathcal{X}] / \mathcal{R}(\delta_{\mathbf{A}}^1) \cong \mathcal{X}/\mathcal{R}(S),
    \end{align*}
    the last two spaces being one-dimensional and spanned by the classes of $1 \otimes e_1$ and of $1 \otimes (e_1 \wedge e_2)$, respectively. For $\mathbf{B}$ we have $\delta_{\mathbf{B}}^0 y = S^* y \otimes e_2$ and $\delta_{\mathbf{B}}^1 (y_1 \otimes e_1 + y_2 \otimes e_2) = - S^* y_1 \otimes (e_1 \wedge e_2)$. Hence $\delta_{\mathbf{B}}^1$ is surjective, $\ker \delta_{\mathbf{B}}^1 = \{ y_1 \otimes e_1 + y_2 \otimes e_2 : y_1 \in \ker S^* \}$, $\mathcal{R}(\delta_{\mathbf{B}}^0) = \{ y \otimes e_2 : y \in \mathcal{Y} \}$, and
    \[
    H^0(\mathbf{B}) = \ker S^*, \qquad H^1(\mathbf{B}) \cong \ker S^*, \qquad H^2(\mathbf{B}) = \{0\},
    \]
    the first two spaces being one-dimensional and spanned by $1$ and by the class of $1 \otimes e_1$, respectively. Thus condition ($i$) of Corollary \ref{cor:necessary} is satisfied, $\delta_{\mathbf{B}}^1$ is surjective, and $H^k(\mathbf{B})$ and $H^{k+1}(\mathbf{A})$ are linearly isomorphic for every $k \in \mathbb{Z}$. In fact, more is true: $\delta_{\mathbf{A}}^0$ is left invertible, with left inverse $x_1 \otimes e_1 + x_2 \otimes e_2 \mapsto S^* x_1$, $\delta_{\mathbf{B}}^1$ is right invertible, with right inverse $w \otimes (e_1 \wedge e_2) \mapsto - Sw \otimes e_1$, and since all the cohomology spaces are at most one-dimensional, the linear isomorphisms $H^k(\mathbf{B}) \cong H^{k+1}(\mathbf{A})$ are topological. Thus the exact analogues of the three hypotheses of Theorem \ref{thm:1D_completion} are satisfied.

    Nevertheless, no coupling tuple $\mathbf{C} = (C_1, C_2)$ makes $\mathbf{M}_{\mathbf{C}}$ a commuting and Taylor invertible pair. Indeed, for the present $\mathbf{A}$ and $\mathbf{B}$ the intertwining relation \eqref{eq:completion_coupling} reads $S C_2 = - C_1 S^*$. Applying it to the constant function $1 \in \ker S^*$ gives $S C_2 1 = 0$, hence $C_2 1 = 0$, since $S$ is injective. Therefore
    \[
    \partial_1^{\mathbf{C}} [1 \otimes e_1] = [\delta_{\mathbf{C}}^1 (1 \otimes e_1)] = [C_1 1 \otimes (e_1 \wedge e_1) + C_2 1 \otimes (e_2 \wedge e_1)] = [- C_2 1 \otimes (e_1 \wedge e_2)] = 0,
    \]
    although $[1 \otimes e_1] \neq 0$ in $H^1(\mathbf{B})$. Thus $\partial_1^{\mathbf{C}}$ is not injective, and $\mathbf{M}_{\mathbf{C}}$ is not Taylor invertible by Theorem \ref{thm:taylor_char}.

    The obstruction is detected by condition ($ii$) of Corollary \ref{cor:necessary}. The map $\varepsilon_1^{\mathbf{B}} \colon H^0(\mathbf{B}) \to H^1(\mathbf{B})$ sends $1$ to the class of $1 \otimes e_1$ and is therefore bijective, whereas $\varepsilon_1^{\mathbf{A}} \colon H^1(\mathbf{A}) \to H^2(\mathbf{A})$ vanishes, because $E_1 (x \otimes e_1) = x \otimes (e_1 \wedge e_1) = 0$ for every $x \in \mathcal{X}$. By Theorem \ref{thm:module}, $H^\bullet(\mathbf{B})$ and $H^{\bullet+1}(\mathbf{A})$ are not isomorphic as graded $\Lambda_2[\mathbf{e}]$-modules, although they are isomorphic, even topologically, as graded vector spaces. Let us also note that the obstruction disappears when the roles of the two variables in $\mathbf{B}$ are interchanged. For $\mathbf{B}' = (S^*, 0)$ the coupling $\mathbf{C} = (I, 0)$ satisfies \eqref{eq:completion_coupling}, and the resulting pair
    \[
    \mathbf{M}_{\mathbf{C}} = \left( \begin{pmatrix} S & I \\ 0 & S^* \end{pmatrix}, \begin{pmatrix} 0 & 0 \\ 0 & 0 \end{pmatrix} \right)
    \]
    is Taylor invertible, since its first entry is invertible, with inverse $(u, v) \mapsto (S^* u - v, Sv + u(0) 1)$.
\end{example}
The following remark compares the multivariable completion problem with its one-dimensional counterpart and isolates the obstruction exhibited by Example \ref{ex:no_completion}.
\begin{remark}\label{rem:obstruction}
    Theorem \ref{thm:completion} should be compared with its one-dimensional counterpart, Theorem \ref{thm:1D_completion}, where $M_C$ admits an invertible completion if and only if $A$ is left invertible, $B$ is right invertible, and $\mathcal{N}(B) \cong \mathcal{X}/\mathcal{R}(A)$. The decisive simplification in the one-dimensional case is that the coupling consists of a single, unconstrained operator $C$: given any Banach-space isomorphism $\psi \colon \mathcal{N}(B) \to \mathcal{X}/\mathcal{R}(A)$, one can construct a bounded operator $C$ so that the connecting homomorphism $\partial_0^{C}$ coincides with $\psi$ (see \cite{Stankovic1D}). Consequently, the criterion is fully intrinsic, depending only on the spaces $\mathcal{N}(B)$ and $\mathcal{X}/\mathcal{R}(A)$ attached to $A$ and $B$.

    In the multivariable setting this freedom is lost. The coupling tuple $\mathbf{C}$ must obey the intertwining relations \eqref{eq:completion_coupling}, and the induced maps on cohomology are constrained to arise from a single Koszul morphism $\delta_{\mathbf{C}} = \sum_{i=1}^n C_i \otimes E_i$ acting in all degrees simultaneously. Example \ref{ex:no_completion} shows that this constraint is effective: the cohomology spaces of $\mathbf{A}$ and $\mathbf{B}$ may be topologically isomorphic in every degree while no Taylor invertible completion exists. Theorem \ref{thm:module} isolates one part of the constraint, namely the compatibility of the connecting homomorphisms with the action of the exterior algebra, and this part already accounts for Example \ref{ex:no_completion}. The remaining part is the realization of a given isomorphism of graded modules by the connecting homomorphisms of one and the same coupling tuple subject to \eqref{eq:completion_coupling}; by Proposition \ref{prop:gauge}, only the class of $\mathbf{C}$ modulo trivial couplings matters here. Equivalently, the question is whether an abstract isomorphism $H^\bullet(\mathbf{B}) \to H^{\bullet+1}(\mathbf{A})$ of graded modules can be promoted to a degree-one cochain morphism of Koszul type between $\mathscr{K}(\mathbf{B})$ and $\mathscr{K}(\mathbf{A})$.
\end{remark}
\begin{remark}\label{rem:finite_dim}
    The completion problem is genuinely infinite-dimensional in nature. Indeed, if $\dim \mathcal{X} < \infty$ and $\dim \mathcal{Y} < \infty$, then a Taylor invertible completion exists if and only if both $\mathbf{A}$ and $\mathbf{B}$ are Taylor invertible. The sufficiency is Corollary \ref{cor:two_out_of_three}($i$) with $\mathbf{C} = \mathbf{0}$. Conversely, if $\mathbf{M}_{\mathbf{C}}$ is Taylor invertible, then $\ker \delta_{\mathbf{A}}^0 = \{0\}$ by Theorem \ref{thm:AB_range_closed}, that is, $0$ is not a joint eigenvalue of $\mathbf{A}$. Since the Taylor spectrum of a commuting tuple of matrices coincides with the set of its joint eigenvalues \cite{ChoTakaguchi84}, $\mathbf{A}$ is Taylor invertible, and Corollary \ref{cor:two_out_of_three}($i$) shows that $\mathbf{B}$ is Taylor invertible as well. The same argument shows that in finite dimensions the necessary conditions of Corollary \ref{cor:necessary} are also sufficient: condition ($i$) forces $\mathbf{A}$ to be Taylor invertible, condition ($ii$) then gives $H^k(\mathbf{B}) = \{0\}$ for all $k$, and $\mathbf{C} = \mathbf{0}$ is a Taylor invertible completion.
\end{remark}

\bigskip
\section*{Declarations}
\vspace{0.5cm}

\noindent{\bf{Funding}}\\
The first author acknowledges the support of the CDP C2EMPI, as well as the French State under the France-2030 programme, the University of Lille, the Initiative of Excellence of the University of Lille, the European Metropolis of Lille for their funding and support of the R-CDP-24-004-C2EMPI project. The second author has been supported by the Ministry of Science, Technological Development and Innovation of the Republic of Serbia [Grant Number: 451-03-34/2026-03/200102].
\vspace{0.5cm}

\noindent{\bf{Availability of data and materials}}\\
\noindent No data were used to support this study.
\vspace{0.5cm}

\noindent{\bf{Competing interests}}\\
\noindent The authors declare that they have no competing interests.
\vspace{0.5cm}

\noindent{\bf{Author contribution}}\\
\noindent The work was a collaborative effort of all authors, who contributed equally to writing the article. All authors have read and approved the final manuscript.
\vspace{0.5cm}

\noindent{\bf{Declaration of generative AI and AI-assisted technologies in the manuscript preparation process}}\\
\noindent During the preparation of this work, the authors used Claude (Anthropic) to identify relevant references, to verify the mathematical arguments, to assist in the derivation and drafting of parts of the results, and for language editing and improvement of the exposition. ChatGPT (OpenAI) was additionally consulted for pre-submission checks. The authors reviewed, verified and edited all generated content and take full responsibility for the content of the published article.

\emergencystretch=1em

\end{document}